\documentclass[11pt,a4paper,oneside]{article}
\usepackage{amssymb,latexsym,amsmath,graphics}

\numberwithin{equation}{section}

\usepackage[left=2.3cm, right=2.3cm, top=2.5cm]{geometry}

\usepackage{mathptmx}  

\usepackage{titlesec} 
\titlespacing*{\section}
{0pt}{10pt}{0pt}
\titleformat{\subsection}[runin]
  {\normalfont\large\bfseries}{\thesubsection}{1em}{}
\titleformat{\subsubsection}[runin]
  {\normalfont\normalsize\bfseries}{\thesubsubsection}{1em}{}

\let\OLDthebibliography\thebibliography
\renewcommand\thebibliography[1]{
  \OLDthebibliography{#1}
  \setlength{\parskip}{0pt}
  \setlength{\itemsep}{0pt plus 0.4ex}
}

\usepackage{graphicx,color}
\usepackage{latexsym}
\usepackage{epsfig}
\usepackage{amssymb}
\usepackage{amstext}
\usepackage{amsgen}
\usepackage{amsxtra}
\usepackage{amsgen}
\usepackage{amsthm}
\usepackage{tocloft} 
\usepackage{authblk}
\usepackage{graphicx}
\usepackage{caption}
\usepackage{subcaption}
\usepackage{enumerate}  
\usepackage{enumitem}
\usepackage{color}
\usepackage{dsfont}
\usepackage{amsmath}
\usepackage{geometry}
\usepackage{marginnote}
\usepackage{bbm}
\usepackage{cite}
\usepackage{authblk}
\usepackage{hyperref} 
\usepackage{mathtools}
\hypersetup{
    pdfpagemode={UseOutlines},
    bookmarksopen,
    pdfstartview={FitH},
    colorlinks,
    linkcolor={blue}, 
    citecolor={blue}, 
    urlcolor={blue} 
}

\newcommand{\pt}{{\widetilde{\psi}}}
\newcommand{\Wt}{{\overline{W}}}
\newcommand{\E}{{\mathbb E}}
\newcommand{\W}{{\mathbb W}}

\newcommand{\R}{{\mathbb R}}
\newcommand{\g}{{\, |\,}}
\newcommand{\Oo}{\mathcal{O}}

\newtheorem {theorem}{Theorem}
\newtheorem {proposition}[theorem]{Proposition}

\newtheorem {lemma}[theorem]{Lemma}
\newtheorem {corollary}[theorem]{Corollary}
\newtheorem{condition}[theorem]{Condition}

\theoremstyle{definition}
\newtheorem{example}[theorem]{Example}
\newtheorem{remark}[theorem]{Remark}

\def\tilde{\widetilde}
\def\hat{\widehat}

\def\B{\mathbb{B}}

\def\N{\mathbb{N}}
\def\P{\mathbb{P}}
\def\R{\mathbb{R}}

\def\W{\mathbb{W}}
\def\1{\mathbbm{1}}

\def\LC{\mathcal{L}}
\def\FC{\mathcal{F}}
\def\NC{\mathcal{N}}
\def\OC{\mathcal{O}}
\def\SC{\mathcal{S}}

\def\supp{\textnormal{supp}}

\newcommand{\eq}[1]{ 
\begin{equation}
\begin{split}
#1
\end{split}
\end{equation}
} 

\newcommand{\eqstar}[1]{
\begin{equation*}
\begin{split}
#1
\end{split}
\end{equation*}
} 

\title{Bernstein--von Mises theorems and uncertainty quantification for linear inverse problems}

\author{Matteo Giordano}
\author{Hanne Kekkonen}
\affil{Centre for Mathematical Sciences, University of Cambridge, \\
Wilberforce Road, Cambridge CB3 0WA, United Kingdom} 
\date{}

\begin{document}
\maketitle

\begin{abstract}
We consider the statistical inverse problem of recovering an unknown function $f$ from a linear measurement corrupted by additive Gaussian white noise. We employ a nonparametric Bayesian approach with standard Gaussian priors, for which the posterior-based reconstruction of $f$ corresponds to a Tikhonov regulariser $\bar f$ with a reproducing kernel Hilbert space norm penalty. We prove a semiparametric Bernstein--von Mises theorem for a large collection of linear functionals of $f$, implying that semiparametric posterior estimation and uncertainty quantification are valid and optimal from a frequentist point of view. The result is applied to study three concrete examples that cover both the mildly and severely ill-posed cases: specifically, an elliptic inverse problem, an elliptic boundary value problem and the heat equation. For the elliptic boundary value problem, we also obtain a nonparametric version of the theorem that entails the convergence of the posterior distribution to a prior-independent infinite-dimensional Gaussian probability measure with minimal covariance. As a consequence, it follows that the Tikhonov regulariser $\bar f$ is an efficient estimator of $f$, and we derive frequentist guarantees for certain credible balls centred at $\bar{f}$.



\end{abstract}

\tableofcontents

%
%
%
%
%

\section{Introduction}
\label{Sec:Intro}

	Inverse problems arise in a variety of scientific disciplines, where the relationship between the quantity of interest and the data collected in an experiment is determined by the physics of the underlying system and can be mathematically modelled. Real world measurements are always discrete and carry statistical noise, which is often most naturally modelled by independent  Gaussian random variables. The observation scheme then gives rise to an inverse regression model of the form
$$
M_i=(Af)_i+ w_i, \quad i=1,\dots,n, \ w_i\overset{\textnormal{iid}}{\sim} \NC(0,1),
$$
where $A$ describes the forward process and $(Af)_i$ is a discrete observation of  the transformed signal.

	The formulation and analysis of the inverse problem is often best done by working with an analogous continuous model. This guarantees, among other things, discretisation invariance that allows to switch consistently between different discretisations \cite{Dashti2016, Lasanen2012, Lassas2009, Stuart2010}. In this paper we consider the case where the forward operator $A:W_1\to W_2$ is linear between separable Hilbert spaces $W_1$ and $W_2$,  and assume the continuous equivalent model (in the sense of \cite{Brown1996b, Reiss2008})
\eq{\label{StatInvProbl}
M_\varepsilon=Af+\varepsilon\W, \quad \varepsilon>0,
}
where $\W$ is a Gaussian white noise process indexed by $W_2$. Note that while $\W$ can be defined by its actions on $W_2$, it almost surely does not take values on it, making the noise in \eqref{StatInvProbl} 'rougher' than the forward signal $Af$.

	We adopt the Bayesian approach to inverse problems \cite{Dashti2016, Stuart2010} and study the performance of nonparametric procedures based on centred Gaussian priors $\Pi$ for $f$. We are interested in what kind of objective guarantees can be achieved for Bayesian inference based on standard Gaussian priors used in practice.  The specification of these priors does not require additional - and often unavailable - information on the forward map $A$, such as its singular value decomposition (SVD). The solution to the statistical inverse problem is the conditional distribution of $f$ given $M_\varepsilon$, whose mean or mode can be used as point estimators. The main appeal of the method is, however, that it  automatically delivers   quantification of uncertainty in the reconstruction, obtained through credible sets, i.e. regions of the parameter space with specified high posterior probability. In many applications this method can be efficiently implemented using modern (possibly infinite-dimensional) MCMC algorithms that allow fast sampling from the posterior distribution \cite{Beskos2017, Kaipio2004a}.

	Our goal is to investigate whether the methodology delivers correct, prior-inde-pendent and possibly optimal inference on the unknown parameter in the small noise limit. These questions can be addressed under the frequentist assumption that  $M_\varepsilon$ is in reality generated through model (\ref{StatInvProbl}) from a fixed true signal $f^\dagger$ (instead of $f$ being randomly drawn from $\Pi$). We then study the asymptotic concentration of the posterior distribution around $f^\dagger$ as $\varepsilon\to0$. The frequentist analysis of nonparametric Bayesian procedures for inverse problems has received increasing interest in the last decade, and several contributions in the linear setting have established consistency results and derived posterior contraction rates; see \cite{Agapiou2013,Agapiou2014,Kekkonen2016,Knapik2011,Knapik2013,Knapik2016,Knapik2018, Ray2013, Vollmer2013} among others. We also mention \cite{Nickl2017,NvdGW2018,Nickl2019} for results for non-linear inverse problems.

	However, determining whether the resulting uncertainty quantification is objectively valid requires finer analysis of the posterior distribution. The central question is: do credible sets have the correct frequentist coverage in the small noise limit? That is, do we have, for some set $C=C(M_\varepsilon)$,
\begin{align}\label{UQvalidation}
\Pi\Big(f\in C \,\Big|\, M_\varepsilon \Big) \approx 1-\alpha \quad\Leftrightarrow\quad
{\mathbb{P}}\Big(f^\dagger\in C\Big)\approx1-\alpha,
\end{align}
with small $\alpha\in(0,1)$ as $\varepsilon\to0$?
The importance of the above questions is not restricted  just to the  Bayesian paradigm. In linear Bayesian inverse problems with Gaussian priors the conditional mean estimator can be shown to coincide with a Tikhonov regulariser $\bar{f}$ arising from a reproducing kernel Hilbert space norm penalty, see \cite{Dashti2013, Helin2015}. Thus, if \eqref{UQvalidation} holds for a credible set $C$ centred at the posterior mean, we can use $C$ as an (asymptotic) frequentist confidence region based on the Tikhonov regulariser $\bar f$.

	Obtaining optimal contraction rates is not enough to answer the above question even in the parametric case. For classical finite-dimensional models the Bernstein--von Mises (BvM) theorem states that, under mild conditions, the posterior distribution is approximated in total variation distance by a normal distribution, centred at the maximum likelihood estimator and with minimal asymptotic variance. This implies that credible sets are asymptotically valid and optimal confidence regions; see, e.g., \cite[Chapter 10]{Vaart1998}. Understanding the  frequentist properties of nonparametric credible sets presents a more delicate matter. It was observed by \cite{Cox1993}, and later in \cite{Freedman1999}, that the theorem may fail to hold even in a simple nonparametric regression model, for which credible balls in $L^2$ can be shown to have null asymptotic coverage.

	One way of tackling the problem is to start by examining the limit behaviour of the one-dimensional marginals $\langle f,\psi\rangle_{W_1}\g M_\varepsilon$ instead of the full posterior. This semiparametric approach was introduced for a direct problem where $A=I$  in \cite{Castillo2013, Castillo2014}, where it is shown that (approximately) in the small noise limit 
\eq{
\label{Eq:Heuristic1}
	\langle f,\psi\rangle_{W_1}|M_\varepsilon \sim\NC(\langle \bar f,\psi\rangle_{W_1}, \varepsilon^2 \mathbb{I}^{-1}(\psi)),
}
for a large collection of test functions $\psi$. Above $\mathbb{I}^{-1}(\psi)$ is the asymptotic minimal variance.  
Note that nonparametric BvM theorems cannot hold in total variation distance like the classical BvM theorem. Instead one has to employ some metric for weak convergence of probability measures. Utilising a Wasserstein-type metric \cite{Castillo2013, Castillo2014} achieve weak convergence of the posterior distribution to a prior-independent infinite-dimensional Gaussian distribution on a large enough function space. 
More recently similar techniques were used 
in the inverse setting \cite{Monard2019}, for the linear X-ray transform problem, obtaining a semiparametric BvM theorem relative to  \textit{smooth} functionals of the unknown, while \cite{Nickl2019}  proved a nonparametric result for a non-linear problem arising in partial differential equations. See also \cite{Nickl2017,NicklRay2019} for further related results. Positive results have also been obtained in \cite{Knapik2011,Leahu2011, Szabo2015}, for priors defined on the SVD basis of the forward operator.

	The first contribution of the present paper is to extend the semiparametric BvM theorem  in \cite{Monard2019} for linear inverse problems of the form \eqref{StatInvProbl},  formulating a general framework that translates the $C^\infty$ smoothness assumption on the test functions $\psi$ into a general 'source-type condition' that depends on the properties of the forward map and of the chosen prior (cf. Theorem \ref{Thm:BvM1}). As a consequence, we then deduce that the plug-in Tikhonov regularisers $\langle \bar f,\psi\rangle_{W_1}$ are consistent and efficient estimators for $\langle f^\dagger,\psi\rangle_{W_1}$, and that credible intervals centred at such estimators have asymptotically correct coverage and optimal width.

	We subsequently employ the general theory to study three concrete examples of interest, where properties of the forward map can be exploited to check the condition for the semiparametric BvM theorem to hold. Specifically, we consider elliptic inverse problems on closed manifolds (Example \ref{ex:psido}), an inverse problem arising from an elliptic boundary value problem (Example \ref{Ex:semiParBvM1}), and the severely ill-posed problem of finding the initial source of the heat equation (Example \ref{Ex:HeatEquation}).   
Similar examples have been considered, e.g, respectively in \cite{Agapiou2013, Kekkonen2016}, in \cite{Gugushvili2018} and in  \cite{Agapiou2014, Knapik2013, Ray2013}.

	Our second contribution is a refinement of the result obtained for the elliptic boundary value problem, for which we further relax the assumption on the test functions to a minimal smoothness requirement that only depends on the degree of ill-posedness (cf.Theorem \ref{thm:semiBvM2}). Adapting the program laid out in \cite{Nickl2017} to the problem at hand, we show that the asymptotic approximation of the marginal distributions holds uniformly across a suitable collection of test functions, leading to the formulation of a nonparametric BvM theorem. This entails the convergence of the posterior distribution to a limiting Gaussian probability measure with minimal covariance in suitable function spaces (cf. Theorem \ref{Mtheo}), and implies frequentist guarantees for the reconstruction and uncertainty quantification relative to the entire function $f$.

	The article is organised as follows:  we introduce the general setting in Section \ref{Sec:BvM}, and state the semiparametric BvM theorem for linear functionals of the unknown in Section \ref{SubSec:BvMtheorem}. In Section \ref{Sec:UQ} we derive the asymptotic normality of $\langle \bar f,\psi\rangle_{W_1}$ and the coverage properties of credible intervals. Section \ref{Sec:Examples} is dedicated to the examples. In Section \ref{Sec:BvMextension} we refine the general theorem to achieve optimal semiparametric result for the elliptic boundary value problem, and obtain the nonparametric BvM theorem. The proofs are postponed to Section \ref{Sec:Proofs} and partly to the Supplement \cite{GK19Supplement} to this article (included below). Finally, Appendix \ref{appendix:infolb} and \ref{Sec:EllipticBVPFacts}  provide some of the background facts used throughout the paper.

	Regarding the notation, we will write $\lesssim$ and $\gtrsim$ for inequalities holding (possibly asymptotically) up to a universal constant. Also, for two real sequences $(a_n)$ and $(b_n),$ we say that $a_n\simeq b_n$ if both $a_n\lesssim b_n$ and $b_n\lesssim a_n$ for all $n$ (large enough). Below, we will denote by $\to^d$ the usual convergence in distribution of a sequence of random variables. The notation $\mu_\varepsilon\to^{\mathcal{L}}\mu$ will be used for the weak convergence of random laws in probability, meaning that for any metric  $d$ for weak convergence of probability measures the real random variables $d(\mu_\varepsilon,\mu)$ converge to zero in probability (see \cite{Dudley2002} for definitions).
		
%
%
%
%
%

\section{General posterior results}
\label{Sec:BvM}

%
%
%

\subsection{A semiparametric Bernstein--von Mises theorem.}
\label{SubSec:BvMtheorem}

	We start by considering general linear inverse problems with minimal assumptions on the forward operator.  We are interested in the nonparametric statistical inverse problem of recovering an unknown function $f$ from a noisy measurement of the form
\begin{align}
\label{eq:observation}
	M_\varepsilon=Af+\varepsilon \W, \quad \varepsilon>0.
\end{align}
The forward operator $A: W_1\to W_2$ is assumed to be linear, bounded and injective between separable Hilbert spaces $W_1 $ and $W_2$ of real valued functions (that can be defined on different sets). The operator $A$ has a well defined adjoint $A^*:W_2\to W_1$ for which $\langle Af, g \rangle_{W_2}=\langle f, A^*g \rangle_{W_1},$ for all $f\in W_1$ and $g\in W_2$. In order to deal with possibly non-smooth unknowns, we define a third space $\Wt$ as a separable Hilbert space for which $A:\Wt\to W_2$ is continuous and $W_1\subset\Wt$ is dense in the norm of $\Wt$. In particular, there exists $c>0$ such that 
\begin{align}\label{eq:SmoothesOfA}
	\|Af\|_{W_2}\leq c\|f\|_\Wt \quad \forall f\in\Wt.
\end{align}

	The above can be thought of as a smoothing property of $A$, in that  the more smoothing the forward operator is, the larger the space $\Wt$ can be chosen. For example, if we assume that $A:L^2(\R^d)\to L^2(\R^d)$ is an elliptic (pseudo-)differential operator smoothing of order $t$, we may choose $\Wt=H^{-t}(\R^d)$, see Section \ref{Subsec:EllipticInvProbl}. Since our general semiparametric result only requires that $f\in\Wt$, this allows dealing with possibly non-smooth unknowns $f\not\in L^2(\R^d)$ as long as $t>0$ (cf. Example \ref{ex:psido} and the following discussion). Note that we can always make the trivial choice $\Wt=W_1$.

	The measurement noise $\W$ is taken to be a centred Gaussian white noise process $(\W(\varphi) : \varphi\in W_2)$ defined on some probability space $(\Omega, \Sigma,\P)$, with covariance 
$
	\E\left( \W(\varphi)\W(\psi) \right) 
	= \langle \varphi,\psi\rangle_{W_2}.
$
Below we often write $\langle \W, \varphi  \rangle_{W_2}$ for the random variable $\W(\varphi)$. 
The noise amplitude is modelled by $\varepsilon>0$. 
Observing data $M_\varepsilon$ then means that we observe a realisation of the Gaussian process $(M_\varepsilon(\varphi)=\langle M_\varepsilon,\varphi \rangle_{W_2} : \varphi\in W_2)$ with marginal distributions $\langle M_\varepsilon,\varphi\rangle_{W_2}\sim \NC(\langle Af,\varphi\rangle_{W_2},\varepsilon^2\|\varphi\|^2_{W_2})$.

	For a fixed $f\in\Wt$, let $P_{f}^M=\mathcal{L}(M_\varepsilon)$ be the (cylindrically defined) law of $M_\varepsilon$. Arguing as in Section 7.4 in \cite{Nickl2019} (see also \cite[Theorem 2.23]{daPrato2014}), we can use the law $P_0^{M}$ of $\varepsilon\W$ as a common dominating measure, and apply the Cameron--Martin theorem \cite[Corollary 2.4.3.]{Bogachev1998} to define the log-likelihood function as
\begin{align}\label{Likelihood}    
	f\mapsto\ell(f,M_\varepsilon)
	=
		\log p_f(M_\varepsilon)
	:=
		\log \frac{dP_{f}^M}{dP_0^M}(M_\varepsilon)  
	= 
		\frac{1}{\varepsilon}\langle M_\varepsilon, 
		Af\rangle_{W_2}-\frac{1}{2\varepsilon^2}\|Af\|_{W_2}^2.
\end{align}

	We consider a Bayesian approach to the problem, assigning $f$ a centred Gaussian prior $\Pi$ on $\Wt$. The reproducing kernel Hilbert space (RKHS) or Cameron-Martin space of $\Pi$ is denoted by $V_\Pi$. Noticing that $\ell(f,M_\varepsilon)$  can be taken to be jointly measurable, we can then use Bayes' theorem to deduce that the posterior distribution of $f\g M_\varepsilon$ arising from observation \eqref{eq:observation} can be written as 	
\begin{align}
\label{eq:BayesFormula}
	\Pi(B\g M_\varepsilon) 
	& = 
		\frac{\int_B p_f(M_\varepsilon)d\Pi(f)}{\int_{\Wt} p_f(M_\varepsilon)d\Pi(f)}
		\quad 
		\text{$B\in\mathcal{B}_{\Wt}$ a Borel set in $\Wt$}.
\end{align}

	In the following we will study the asymptotic behaviour of $\Pi(\cdot \g M_\varepsilon)$ in the small noise limit $\varepsilon\to0$, under the assumption that the measurement is generated from a fixed true unknown $f^\dagger\in\Wt$. In order to do so, we assume that the prior satisfy a standard concentration function condition.

\begin{condition}\label{cond:concentration}

	Let $\Pi$ be a centred Gaussian Borel probability measure on the separable Hilbert space $\Wt$ for which \eqref{eq:SmoothesOfA} holds,
and let $V_\Pi$ be the RKHS of $\Pi$. Define the concentration function of $\Pi$ with a fixed $f^\dagger\in\Wt$ as
\begin{align}\label{eq:concentration}
	\phi_{\Pi,f^\dagger}(\delta)
	=
		\inf_{g\in V_\Pi,\ \|g-f^\dagger\|_{\Wt}\leq\delta}\frac{\|g\|_{V_\Pi}^2}{2}
		-\log\Pi(f : \|f\|_\Wt\leq\delta), \quad \delta>0.
\end{align}
Given $\Pi$ and $f^\dagger\in \Wt$ assume that there exists a sequence $\delta_\varepsilon\to0$, with ${\delta_\varepsilon}/{\varepsilon}\to\infty$ as $\varepsilon\to0$, such that
\begin{align}
\label{eq:concentrationAs}  
	\phi_{\Pi,f^\dagger}\left(\frac{\delta_\varepsilon}{2c}\right)
	\leq
		 \left(\frac{\delta_\varepsilon}{\varepsilon}\right)^2.
\end{align}
\end{condition}

	The above condition characterises the asymptotics of the small ball probabilities, and guarantees that the prior puts sufficient mass around the truth: in particular $\Pi(f : \|f-f^\dagger\|_{\Wt}\leq\delta_\varepsilon)>e^{-\frac{1}{2}(\delta_\varepsilon/\varepsilon)^2}$ as $\delta_\varepsilon\to0$ (cf. the proof of Lemma \ref{App:ProofLemAe}). Analogous conditions underpin many results in Bayesian asymptotics, and play a fundamental role in the theory of posterior contraction rates, see e.g. \cite{Ghosal2017,Gine2016,Vaart2008}. The concentration functions of Gaussian priors are generally well understood, and explicit forms for the sequences $\delta_\varepsilon$ can readily be computed for many standard choices of practical interest, such as the commonly used Mat\'ern process priors (see Section \ref{Sec:Examples}).

	Next we formulate a semiparametric Bernstein--von Mises theorem in the above general linear inverse problems setting.
  
\begin{theorem}\label{Thm:BvM1}

Let $P_{f^\dagger}^M$ be the law of $M_\varepsilon$ generated by  \eqref{eq:observation} with $f=f^\dagger\in\Wt$ , where $\Wt$ is a separable Hilbert space for which \eqref{eq:SmoothesOfA} holds.  
 We assume a centred Gaussian prior $\Pi$ that satisfies Condition \ref{cond:concentration} for a fixed $f^\dagger\in\Wt$ and denote its RKHS by $V_\Pi$. 
	 Consider a test function $\psi\in W_1$ such that $|\langle\psi,\varphi\rangle_{W_1}|\lesssim \|\varphi\|_\Wt$, for all $\varphi\in W_1$,
and suppose that $\psi=-A^*A\pt$ for some $\pt \in V_\Pi$. Then,
\begin{align}\label{lawConv}
	\mathcal{L}\left(\varepsilon^{-1}(\langle f,\psi\rangle_{W_1}-\widehat{\Psi}) \g 	 M_\varepsilon\right) \to ^\mathcal{L} \NC(0,\|A\pt\|_{W_2}^2)
\end{align}
in $P_{f^\dagger}^M$-probability as $\varepsilon\to0$, where
\begin{align}\label{HatPsi}
	\widehat{\Psi}
	=
		\langle f^\dagger,\psi\rangle_{W_1}-\varepsilon\langle A\pt, \W\rangle_{W_2}.
\end{align}
\end{theorem}

	The next corollary states that we can replace the centring $\widehat{\Psi}$ by a linear functional of the conditional mean. This implies that the posterior distribution of the functionals
are asymptotically approximated by a normal distribution centred at the conditional mean and with asymptotic minimal variance (see Remark \ref{TikEfficiency} below). The proof of Corollary \ref{Thm:BvM2} can be adapted from the proof of Theorem 2.7 in \cite{Monard2019} and is therefore omitted (see also Step V in the Supplement).

\begin{corollary}\label{Thm:BvM2}

	Let $\bar f=\E^{\Pi}[f|M_\varepsilon]$ be the mean of the posterior $\Pi(\cdot|M_\varepsilon)$. Then, for every $\psi\in W_1$ satisfying the conditions in Theorem \ref{Thm:BvM1}, we have
\eq{\label{meanConv}
	\frac{1}{\varepsilon}(\langle \bar f,\psi\rangle_{W_1}-\hat\Psi) \to0,
}
in $P^M_{f^\dagger}-$probability as $\varepsilon\to0$. As a consequence, we can replace $\hat\Psi$ with $\langle \bar f,\psi\rangle_{W_1}$ in Theorem  \ref{Thm:BvM1}.
\end{corollary}

	Note that, since $W_1\subset\Wt$ is dense and $L_\psi(\cdot)=\langle\psi,\cdot\rangle_{W_1}$ is assumed to be a bounded linear operator (and hence uniformly continuous), we can extend $L_\psi$ continuously to $\Wt$.  The condition on the test functions requires that $\psi$ is in the range of the 'Fisher information operator' $A^*A$ acting upon the RKHS of $\Pi$. This can normally be translated into suitable smoothness assumptions on $\psi$, see Section \ref{Sec:Examples} for  examples. The requirement resembles certain source conditions often used in inverse problems \cite{Engl1996a, Lu2013, Schuster2012}.
The main conceptual difference is that instead of requiring extra smoothness for the unknown $f^\dagger$ to attain convergence in a predefined space, we allow $f^\dagger$ to be non-smooth and impose constraints on the test functions in order to achieve convergence.

%
%
%

\subsection{Efficiency and uncertainty quantification for Tikhonov regularisers.}
\label{Sec:UQ}

	Since the forward operator $A$  is assumed to be linear, the posterior distribution $\Pi(\cdot|M_\varepsilon)$ is Gaussian.
 It follows that the conditional mean $\bar f=\bar f(M_\varepsilon):=\E^\Pi[f|M_\varepsilon]$
coincides with the maximum a posterior (MAP) estimator, and using  Corollary 3.10 in \cite{Dashti2013} (under appropriate conditions on $A$) the latter can be seen to be a  Tikhonov-type regulariser found by minimising the following Onsager-Machlup functional
$$
Q(f)=-\frac{1}{\varepsilon^2}\langle M_\varepsilon,Af\rangle_{W_2}+\frac{1}{2\varepsilon^2}\|Af\|^2_{W_2}+\frac{1}{2}\|f\|^2_{V_\Pi}.
$$

	Using Theorem \ref{Thm:BvM1} and Corollary \ref{Thm:BvM2} we can derive the asymptotic distribution of the plug-in estimators $\langle \bar f,\psi\rangle_{W_1}$. 

\begin{remark}[Minimax optimality of the plug-in Tikhonov regulariser]\label{TikEfficiency}
Corollary \ref{Thm:BvM2} implies that 
\eq{
\label{meanConv2}
	\frac{1}{\varepsilon}\langle \bar f-f^\dagger,\psi\rangle_{W_1}
	\to^{d} Z\sim\NC(0,\|A\tilde\psi\|^2_{W_2})
}
in $P^M_{f^\dagger}-$probability as $\varepsilon\to0$. The above random variable $Z$  identifies the asymptotic minimal variance (in the minimax sense) in estimating $\langle f^\dagger,\psi \rangle_{W_1}$ from model (\ref{eq:observation}), in that 
\eq{
\label{infolb1}
	\liminf_{\varepsilon\to0}\inf_{T}\sup_{f\in B_\varepsilon}\varepsilon^{-2}
	E_{f^\dagger}^M(\langle f^\dagger,\psi\rangle_{W_1}-T)^2
\ge
\|A\tilde\psi\|^2_{W_2},
}
the infimum being over all estimators $T=T (M_\varepsilon,\psi)$ of $\langle f^\dagger,\psi\rangle_{W_1}$ based on observing $M_\varepsilon$ in (\ref{eq:observation}) with $f=f^\dagger$, and the supremum is taken over balls $B_\varepsilon$ in $\overline W$ centred at $f^\dagger$ and with radius $\varepsilon>0$; see Appendix \ref{appendix:infolb}.

	We notice that (\ref{meanConv2}) implies the convergence of all moments (see Step V in the Supplement). Consequently, for all $\psi\in W_1$ fulfilling the conditions of Theorem \ref{Thm:BvM1}, the {plug-in} Tikhonov regulariser $\langle \bar f,\psi\rangle_{W_1}$ attains the lower bound in (\ref{infolb1}), and hence is an asymptotic minimax estimator of $\langle f^\dagger,\psi\rangle_{W_1}$.

\end{remark}

	Besides the question of efficiency, the most relevant consequence of Theorem \ref{Thm:BvM1} is that credible intervals built around the estimators $\langle \bar f,\psi\rangle_{W_1}$ are asymptotically valid frequentist confidence intervals with optimal diameter. Specifically, for $\psi$ as above, consider a credible interval for $\langle \bar f,\psi\rangle_{W_1}$ of the form
\eq{\label{Ceps}
	C_\varepsilon=\{x\in\R : |\langle \bar f,\psi\rangle_{W_1}-x|\le R_\varepsilon\},
}
with $R_\varepsilon= R_\varepsilon(\alpha,M_\varepsilon)$ chosen so that
$$
	\Pi(\langle f,\psi\rangle_{W_1}\in C_\varepsilon \g M_\varepsilon)
	=
	1-\alpha,
	\quad \alpha\in(0,1).
$$  
Then it follows that $C_\varepsilon$ has the correct asymptotic coverage and that its diameter shrinks at the optimal rate $\varepsilon$. The proof of the following corollary can be found in the Supplement.

\begin{corollary}\label{Cor:SemiUQ}
 
 	Let $\psi\in W_1$ satisfy the conditions in Theorem \ref{Thm:BvM1}, and let $C_\varepsilon$ be as in (\ref{Ceps}). Then, as $\varepsilon\to0,$
\begin{align*}
	P^M_{f^\dagger}(\langle f^\dagger,\psi\rangle_{W_1}\in C_\varepsilon) 
	& \to 1-\alpha\quad\textnormal{and}\\
	\varepsilon^{-1} R_\varepsilon & \to^{P^M_{f^\dagger}}\Phi^{-1}(1-\alpha).
\end{align*}
Here $\Phi(t)=\Pr(|Z|\le t) $ and $Z\sim \NC(0,\|A\pt\|^2_{W_2})$.
\end{corollary}

	Note that although an explicit formulation of $C_\varepsilon$ would require the computation of the quantiles of the posterior distribution of $\langle f,\psi\rangle_{W_1}|M_\varepsilon$, these type of credible intervals can often in practice be implemented by numerically approximating the radius $R_\varepsilon$ with a posterior sampling method. See, e.g., \cite{Kaipio2004a}, or Section 2.2 in \cite{Monard2019}.

	For the inferential problem for elliptic partial differential equations studied in Section \ref{sec:ellipticBVP}, Remark \ref{Remark:UQapplications} below will extend the conclusions of Corollary \ref{Cor:SemiUQ} to entire credible balls in suitable function spaces centred at $\bar f$.

%
%
%
%
%

\section{Examples}
\label{Sec:Examples}

	In this section we consider examples of linear inverse problems fitting in the framework of Section \ref{Sec:BvM}, studying the conditions under which the semiparametric Bernstein--von Mises phenomenon occurs in such instances. We first need to introduce some notation on Sobolev spaces (see \cite{Lions1972,McLean2000} for background).

	The Sobolev space on $\R^d$ of order $s\in\R$ is defined as 
\eq{
\label{SobSpaceRd}  
	H^s(\R^d)=\{u\in\SC'(\R^d)\ :\ (1+|\cdot|^2)^{s/2}\FC u\in L^2(\R^d)\},
}
where $\SC'(\R^d)$ is the space of tempered distributions on $\R^d$ and $\FC$ is the Fourier transform. For $\OC\subset\R^d$ a non-empty, open and bounded set  with smooth boundary $\partial\OC$ (a smooth domain), Sobolev spaces on $\OC$ can be defined via the restriction operator $|_{\OC}$ as
\eq{\label{sob2}
	H^s(\OC)
	=
		\{u=U|_{\OC}, \ U\in H^s(\R^d)\}, 
		\quad \|u\|_{H^s(\OC)}=\inf_{U\in H^s(\R^d),\ U|_{\OC}=u}\|U\|_{H^s(\R^d)}.
}

	To correctly address issues relative to the behaviour of functions near $\partial\OC$, we will need to consider certain subspaces of $H^s(\OC)$. We denote the set of functions in $H^s(\OC)$ that are compactly supported in $\OC$ by
$
H^s_c(\OC),
$
and for any fixed compact subset $K\subset\OC$, we write
$
	H^s_K(\OC):=\{u\in H^s(\OC), \ \supp(u)\subseteq K\}.
$
Finally, for all $s>1/2$, let
$
	H^s_0(\OC)
$
be the usual subspace of $H^s(\OC)$ of functions with null trace on $\partial\OC$. Below we will often suppress the dependence on the underlying domain denoting $H^s=H^s(\OC)$.


%
%
%

\subsection{Elliptic inverse problems.}
\label{Subsec:EllipticInvProbl}

	We start with a basic example to demonstrate how Theorem \ref{Thm:BvM1} can be applied when $A$ is assumed to be a smoothing elliptic pseudo-differential operator and $\Oo$ a closed manifold (see \cite{Hormander1994, Shubin1987} for general theory on pseudo-differential operators). The previous definitions of Sobolev spaces can straightforwardly be adapted to this setting, see e.g. \cite[Chapter I.7]{Shubin1987}. The absence of a boundary and the properties of the forward map allow for a clean exposition of the results. In the Section \ref{sec:ellipticBVP} we will instead assume that $\OC$ is a smooth domain in $\R^d$, and take $A$ to be the solution operator associated with an elliptic boundary value problem. We then have to refine the results to take into account some subtleties of the behaviour of functions near the boundary. 

\begin{example}\label{ex:psido}

	Let $\Oo$ be a closed $d$-dimensional manifold and $A:L^2(\Oo)\to L^2(\Oo)$ an injective and elliptic pseudo-differential operator smoothing of order $t$, that is, $A:H^s(\Oo)\to H^{s+t}(\Oo)$ with all $s\in\R$ \cite[Section I.5.]{Shubin1987}. We can then choose $\Wt=H^{-t}(\Oo)$. 
	
	Let $P^M_{f^\dagger}$ be the law of $M_\varepsilon$ generated by \eqref{eq:observation} with $f=f^\dagger\in H^\alpha(\Oo)$, $\alpha>-t$. We assume a centred Gaussian prior $\Pi$ with RKHS $V_\Pi=H^r(\Oo)$, where $r\geq\max\{0,d_0-t\}$ and $d_0>d/2$. This guarantees that $f\in H^{r-d_0}(\Oo)\subset H^{-t}(\Oo)=\Wt$ almost surely. For example, we can take $\Pi= \NC(0,C_f)$, where $C_f$ is a self-adjoint, injective and elliptic covariance operator smoothing of order $2r$ \cite{Agapiou2013, Kekkonen2016}.
Another example is to assume $\Pi$ to be the law of the Mat\'ern process of smoothness $r-d/2$ (see Example 11.8 in \cite{Ghosal2017} for details), namely the centred Gaussian process $(M(x) : x\in\Oo)$ with covariance kernel 
$$
	K(x,y)
	=
		\int_{\R^d}e^{-i\langle x-y,\xi\rangle_{\R^d}} \mu(d\xi),
	\quad
	 \mu(d \xi)=(1+|\xi |^2)^{-r}d\xi.
$$ 
 	
Since $A$  is elliptic and $\OC$ is a closed manifold, $A^*A$ has a well defined inverse $(A^*A)^{-1}:H^s(\Oo)\to H^{s-2t}(\Oo)$,  $s\in\R$, see e.g. \cite{Kekkonen2014}. We can then take $\psi\in H^{r+2t}(\Oo)$, which guarantees $\pt=-(A^*A)^{-1}\psi\in H^r(\Oo)=V_\Pi$ and $|\langle \psi, \varphi\rangle_{L^2}|\leq C\|\varphi\|_{H^{-t}}$, for all $\varphi\in L^2(\Oo)$.

Denote by $\bar f=\E^\Pi[f|M_\varepsilon]$ the mean of the posterior distribution $\Pi(\cdot|M_\varepsilon)$ arising from observing (\ref{eq:observation}). Then, for all test functions $\psi\in H^{r+2t}(\Oo)$, the following convergence occurs in $P^M_{f^\dagger}$-probability as $\varepsilon\to0$
\begin{align*}
	\mathcal{L}\left(\varepsilon^{-1}\langle f-\bar{f},\psi\rangle_{L^2} \g M_\varepsilon\right)
	 \to ^\mathcal{L} \NC\left(0,\|A(A^*A)^{-1}\psi\|_{L^2}^2\right).
\end{align*}

\end{example}

	Note that if $t>\frac{d}{2}-1$ we can allow unknowns of bounded variation $f^\dagger\in BV(\OC)$, since $BV(\OC) \subset H^{\alpha}(\Oo)$ when $\alpha\leq1-\frac{d}{2}$. Functions of bounded variation are widely used e.g. in image analysis due to their ability to deal with discontinuities. One standard example is total variation denoising \cite{Burger2004, Rudin1992}. 

\begin{remark}\label{Rem:EllipticSpeed}
Let $\Pi$ and $f^\dagger$ be as above. Then, as $\delta\to0$
\eq{
\label{eq:speed}
	\phi_{\Pi,f^\dagger}(\delta)
	\lesssim
		\delta^{-\frac{2\max\{0,r-\alpha\}}{t+\alpha}} + \delta ^{-\frac{d}{r+t-d/2}},
}
so that the concentration condition 
$\phi_{\Pi,f^\dagger}(\delta_\varepsilon)\lesssim(\delta_\varepsilon/\varepsilon)^2$
is satisfied by taking
\begin{align*}
	\delta_\varepsilon
	\simeq
		\max\left\{\varepsilon^\frac{t+\alpha}{t+r}, \varepsilon^\frac{t+r-d/2}{t+r}\right\}.
\end{align*}
\end{remark}

	The proof of Remark \ref{Rem:EllipticSpeed} is omitted since it is a simplified version of the proof of Remark \ref{Rem:ConcFunAsymptForL} where, $\Oo$ being a closed manifold, one does not need to address the technicalities arising at the boundary.

%
%
%

\subsection{An elliptic boundary value problem.}
\label{sec:ellipticBVP}

	Let $\OC\subset\R^d$ be a non-empty, open and bounded set with smooth boundary $\partial\OC$. We consider the problem of recovering the unknown source $f\in L^2=L^2(\OC)$ in the elliptic boundary value problem (BVP)
\eq{
\label{BVP}
	\begin{cases}
	Lu=f& \textrm{on}\ \OC\\
	u=0 & \textrm{on}\ \partial\OC
	\end{cases}
}  
from noisy observations of the solution $u$ corrupted by additive Gaussian white noise in $L^2$. We take $L$ to be the following partial differential operator in divergence form:
\eq{
\label{ellipticop}
	Lu=-\sum_{i,j=1}^d\frac{\partial}{\partial x_j}\left(a_{ij}\frac{\partial u}{\partial x_i}\right),
}
for known $a_{ij}\in C^{\infty}(\overline{\OC})$, with $a_{ij}=a_{ji}$. The problem represents an 'elliptic counterpart' of the transport PDE arising in \cite{Monard2019}.

	Assuming that $L$ is uniformly elliptic (see Appendix \ref{Sec:EllipticBVPFacts}), it follows that for each $f\in H^s, \ s\ge0$, there exists a unique weak solution $L^{-1}f\in H^{s+2}_0$ to \eqref{BVP}. In particular,  $L^{-1}:H^s\to H^{s+2}_0$ defines a bounded isomorphism,  self-adjoint with respect to $\langle\cdot,\cdot\rangle_{L^2}$, and for all $s\ge0$ we also have the dual estimates
\begin{align}
\label{DualSobEstim}
	\|L^{-1}f\|_{ (H^s)^*}
	=
		\sup_{u\in H^s,\ \|u\|_{H^s}\le1}|\langle L^{-1}f,u\rangle_{L^2}|
	\le 
		c_s\|f\|_{(H_0^{s+2})^*}
	\quad 
		\text{for some $c_s>0$.}
\end{align}

	Rephrasing in the notation of Section \ref{Sec:BvM}, we consider the observation
\eq{
\label{invprobl1}
M_\varepsilon = L^{-1}f+\varepsilon\W, \quad \varepsilon>0,
}
where $\W$ is Gaussian white noise in $L^2$. For $W_1=W_2=L^2$, the dual estimate (\ref{DualSobEstim}) implies that we can take $\overline W=(H^2_0)^*$.

	We assume that $f\sim\Pi$, where $\Pi$ is a centred Gaussian Borel probability measure on $L^2$ with RKHS $V_\Pi=H^r$, for some $r> d/2$. For example, we can take $\Pi$ to be the law of the Mat\'ern process of smoothness $r-d/2$ introduced in the previous example.

	For  $f^\dagger\in H_c^\alpha$, with some $\alpha\ge 0$, we show that the semiparametric BvM phenomenon occurs under appropriate smoothness conditions on the test functions $\psi$. In particular, assuming that $\psi\in H^{r+4}_c$ automatically verifies the requirements of Theorem \ref{Thm:BvM1}, since taking $\tilde\psi=-L(L\psi)$ implies $\tilde\psi\in V_\Pi = H^r$ and $\psi=-L^{-1}L^{-1}\tilde\psi$, as $\supp(L\psi)\subseteq \supp(\psi)\subsetneq \OC$. The proof of the following proposition can be found in Section \ref{Sec:ProofExample}.
	
\begin{proposition}\label{Ex:semiParBvM1}

	Let $\Pi$ be a Gaussian Borel probability measure on $L^2(\OC)$ with RKHS $V_\Pi=H^r(\OC), \ r>d/2$. Assume that $f^\dagger\in H_c^\alpha(\OC), \ \alpha\ge 0$, and let $P^M_{f^\dagger}$ be the law of $M_\varepsilon$ generated by (\ref{invprobl1}) with $f=f^\dagger$. Let $\bar f=\E^\Pi[f|M_\varepsilon]$ be the mean of the posterior distribution $\Pi(\cdot|M_\varepsilon)$ arising from observing \eqref{invprobl1}. Then, for all $\psi\in H^{r+4}_c(\OC)$, we have 
\eq{\label{finalconv1}  
	\LC(\varepsilon^{-1}\langle f-\bar f,\psi\rangle_{L^2}|M_\varepsilon)
	\to^\LC\NC(0,\|L\psi\|^2_{L^2})
}
in $P^M_{f^\dagger}$-probability as $\varepsilon\to0$.

\end{proposition}

\begin{remark}\label{Rem:ConcFunAsymptForL}
	
Let $\Pi$ and $f^\dagger$ be as above. In the proof of Proposition \ref{Ex:semiParBvM1} we show that as $\delta\to0$
\eq{\label{ConcFunAsymptForL}
	\phi_{\Pi,f^\dagger}(\delta)
	\lesssim
		\delta^{-\frac{2\max\{0,r-\alpha\}}{2+\alpha}} + \delta ^{-\frac{d}{r+2-d/2}},
}
so that the concentration condition 
$\phi_{\Pi,f^\dagger}(\delta_\varepsilon)\lesssim(\delta_\varepsilon/\varepsilon)^2$
is satisfied by taking
\eq{\label{DeltaEpsForL}
	\delta_\varepsilon
	\simeq
		\max\{\varepsilon^\frac{2+\alpha}{2+r},\varepsilon^\frac{2+r-d/2}{2+r}\}.
}

\end{remark}

\subsection{Boundary value problem for the heat equation.}
\label{Subsec:HeatEquation}

	We will conclude this section by applying the general framework studied in Section \ref{Sec:BvM} to the severely ill-posed problem of finding the initial source of  the heat equation. Contraction rates for	similar inverse problems have been studied in \cite{Agapiou2014, Knapik2013, Ray2013}.

\begin{example}\label{Ex:HeatEquation}

	Let $\Oo\subset\R^d$ be an open bounded set with $C^\infty$ boundary $\partial\Oo$. We consider the boundary value problem for the heat equation
\begin{equation*}
\begin{cases}
	u_t-\Delta u =0 & \text{on $\Oo\times\R^+$}\\
	u  = 0 & \text{on $\partial\Oo\times\R^+$}\\
	u(\cdot,0)  = f &  \text{on $\Oo$}.
\end{cases}
\end{equation*}
The inverse problem is to recover the initial heat source $f\in L^2$ from a noisy observation of the solution $u$ at time $T$, corrupted by additive Gaussian white noise on $L^2$. The solution to the boundary value problem is given by 
\begin{align*} 
	u(x,T)=Af(x)
	=
		 \sum_{j=1}^\infty \langle f,\varphi_j\rangle_{L^2} e^{-\lambda_jT} \varphi_j(x), 
	\quad x\in\OC,
\end{align*}  
where $-\Delta\varphi_j=\lambda_j\varphi_j$, and $\{\varphi_j\}_{j=1}^\infty$ forms an orthonormal basis of $L^2$.
If we order the eigenvalues to be increasing, that is, $\lambda_1\leq\lambda_2\leq\dots,$ then Weyl's law yields that $\lambda_j\simeq j^{2/d}$ (e.g., \cite[Theorem 8.16]{Roe1998}). Thus, the singular values of the compact forward operator $A$ decay exponentially to zero, meaning that the recovery of the initial condition of the heat equation is a severely ill-posed inverse problem.

	Assume that $f\sim\Pi$, where $\Pi$ is a centred Gaussian Borel probability measure on $L^2$ with RKHS $V_\Pi=H^r$, $r> d/2$. Let $\psi\in L^2$ be of the form 
\begin{align}
\label{Eq:AnalyticAssumption}
	\psi 
	=
	 	-A^*A\pt
	=
		-\sum_{j=1}^\infty\langle\pt,\varphi_j\rangle_{L^2}e^{-2\lambda_jT}\varphi_j,
\end{align}
for some $\pt\in H^r$. Then
\begin{align*}
	|\langle\psi,\phi\rangle_{L^2}| 
	&=
		 |\langle\psi,\sum_{j=1}^\infty\langle\phi,\varphi_j\rangle_{L^2}\varphi_j\rangle_{L^2}|		\\
	& \leq
		\sum_{j=1}^\infty|\langle\phi,\varphi_j\rangle_{L^2}|\left|\left\langle
		\sum_{i=1}^\infty\langle\pt,\varphi_i\rangle_{L^2}e^{-2\lambda_iT}\varphi_i,			\varphi_j\right\rangle_{L^2}\right|\\
	& \leq
	 	\sum_{j=1}^\infty|\langle\phi,\varphi_j\rangle_{L^2}|
		|\langle\pt,\varphi_j\rangle_{L^2}|e^{-2\lambda_jT}\\
	& \leq 
		C\|\phi\|_{H^{-t}}, 
\end{align*}
for all $t\ge0$, verifying the condition of Theorem \ref{Thm:BvM1}. 
Hence, for $f^\dagger\in L^2$ and $\psi$ as above, we get the following convergence in $P^M_{f^\dagger}$-probability as $\varepsilon\to0$
\begin{align*}
\mathcal{L}\left(\varepsilon^{-1}\langle f-\bar{f},\psi\rangle_{L^2} \g M_\varepsilon\right) \to ^\mathcal{L} \NC\left(0,\|A\pt\|_{L^2}^2\right).
\end{align*}

	Note that the contraction rate $\varepsilon$ entailed by the semiparametric BvM theorem  is a very strong requirement for severely ill-posed inverse problems, usually characterised by logarithmic rates even for smooth functionals \cite{Knapik2013}. To achieve the rate $\varepsilon$, we then need to assume the analytic-type condition \eqref{Eq:AnalyticAssumption} on the test function $\psi$, which  reflects the natural condition of $\psi$ being in the range of $A^*$ which is necessary for efficient semiparametric estimation  \cite[Theorem 25.32]{Vaart1998}.

\end{example}

%
%
%
%
%

\section{A nonparametric Bernstein--von Mises theorem for elliptic boundary value problems}
\label{Sec:BvMextension}

	In this section we continue the investigation of the BvM phenomenon in the setting of the elliptic BVP studied in Section \ref{sec:ellipticBVP}. We develop Proposition \ref{Ex:semiParBvM1} along two related directions: first, we extend the class of test functions $\psi$ for which the convergence (\ref{finalconv1}) occurs, identifying a natural lower limit for the smoothness of $\psi$ that only depends on the level of ill-posedness of the inverse problem. Secondly, combining the result with the program laid out in \cite{Nickl2017}, we derive a nonparametric BvM theorem that entails the  weak convergence, in a suitable function space, of the centred and scaled posterior to a prior-independent infinite-dimensional Gaussian probability measure whose covariance function attains the information lower bound. From the latter result we then obtain frequentist guarantees for uncertainty quantification in the reconstruction of the entire function $f$.

	We briefly recall that, for unknown $f\in L^2=L^2(\OC)$, we consider observations
$
	M_\varepsilon = L^{-1}f+\varepsilon\W,\ \varepsilon>0,
$
where $L^{-1}$ is the solution map associated with the BVP (\ref{BVP}) (see Section \ref{sec:ellipticBVP} for details) and $\W$ is a Gaussian white noise in $L^2$. We assign $f$ a centred Gaussian prior in $L^2$ with RKHS $H^r$, $r> d/2$, and assume that the observation $M_\varepsilon$ is generated from a fixed $f^\dagger\in H^\alpha_c$ with some $\alpha>0$. For the results in this section we assume an undersmoothing prior. That is, we consider the case $r-d/2\le\alpha$. The proofs can be found in Section \ref{Sec:Proofs}. 

\begin{theorem}\label{thm:semiBvM2}

	Let $\Pi$ be a Gaussian Borel probability measure on $L^2(\OC)$ with RKHS $V_\Pi=H^r(\OC), \ r>d/2$. Assume that $f^\dagger\in H_c^\alpha(\OC), \ \alpha\ge r-d/2$, and let $P^M_{f^\dagger}$ be the law of $M_\varepsilon$ generated by (\ref{invprobl1}) with $f=f^\dagger$. Let $\bar f=\E^\Pi[f|M_\varepsilon]$ be the mean of the posterior distribution $\Pi(\cdot|M_\varepsilon)$ arising from observing \eqref{invprobl1}. Then, for all $\beta>2+d/2$, and any $\psi\in H^\beta_c(\OC)$, we have
\eq{
\label{finalconv11}  
	\LC(\varepsilon^{-1}\langle f-\bar f,\psi\rangle_{L^2}|M_\varepsilon)
	\to^\LC \NC(0,\|L\psi\|^2_{L^2})
}
in $P^M_{f^\dagger}$-probability as $\varepsilon\to0$.

\end{theorem}

	Assuming that $\beta>2+d$, we will strengthen the above result to a nonparametric Bernstein--von Mises theorem in the dual spaces $(H_K^\beta)^*$, for any compact set $K\subset \OC$.  In particular, we notice that the Gaussian laws in the right hand side of \eqref{finalconv11} identify the one-dimensional marginal distributions of a nonparametric Gaussian probability measure $\mu$, induced via Kolmogorov's extension (see, e.g., \cite[Section 12.1]{Dudley2002}) on the cylindrical $\sigma$-field of $\R^{H^\beta_K}$ by the centred Gaussian process
\eq{\label{X}
	X=(X(\psi) : \psi\in H^\beta_K),
	\quad 
	E[X(\psi)X(\psi')]=\langle L\psi,L\psi'\rangle_{L^2}.
}

	In fact, in the view of the efficiency considerations in Appendix \ref{appendix:infolb}, $\mu$ represents the 'canonical' asymptotic distribution for the problem of inferring $f$ in model (\ref{invprobl1}), as its covariance function is minimal in the information theoretic sense of Remark \ref{TikEfficiency}. In the following lemma  we derive the values of $\beta$ for which $\mu$ is a tight  Borel probability measure on $(H^\beta_K)^*$, a necessary condition for any sequence of laws on such spaces to weakly converge to $\mu$. 
The proof adapts the argument in the proof of Proposition 6 in \cite{Nickl2017}, and is included in the Supplement.

\begin{lemma}\label{lemma:tightness}

	Fix any compact set $K\subset\OC$. Let $X$ be as in (\ref{X}), and let  $\mu$ be the law of $X$ on the cylindrical $\sigma$-field of $\R^{H^\beta_K}$. Then,
\begin{enumerate}
\item[1.] for all $\beta>2+d/2$, $\mu$ is a tight Gaussian Borel probability measure on $(H^\beta_K)^*$;

\item[2.] for $\beta<2+d/2$, we have
$$
\mu\big(x:\|x\|_{(H^\beta_K)^*}<\infty\big)=0;
$$
\item[3.] for $\beta=2+d/2$, $\mu$ is not tight on $(H^\beta_K)^*$. 
\end{enumerate}

\end{lemma}

	Similarly, the stochastic process obtained by collecting the random variables in the left hand side of \eqref{finalconv11},
\eq{\label{postsp}
X_\varepsilon=\left(\varepsilon^{-1}\langle f-\bar f,\psi\rangle_{L^2}|M_\varepsilon : \psi\in H^\beta_K\right), \quad \varepsilon>0,
}
can also be shown to induce a tight Borel probability measure on $(H^\beta_K)^*$ when $\beta>2+d/2$ (see Step IV in Section \ref{sec:MtheoProof}). We will interpret the law of $X_\varepsilon$ as the nonparametric centred and scaled posterior distribution arising from observing \eqref{invprobl1}, denoted by
\eq{\label{centrScalPost}
\LC(\varepsilon^{-1}(f-\bar f)|M_\varepsilon)=\LC(X_\varepsilon).
}

	Theorem \ref{thm:semiBvM2} implies the convergence of the finite-dimensional distributions of the stochastic process $X_\varepsilon$ to those of $X$ (cf. Lemma \ref{lemma:fddconv}), and by showing that \eqref{finalconv11} holds \textit{uniformly} across the set of test functions, we then deduce the weak convergence of the respective induced laws on $(H^\beta_K)^*$. As mentioned in the introduction, nonparametric BvM theorems cannot hold in total variation distance like the classical BvM theorem. Instead we use a Wasserstein-type metric for weak convergence of probability measures. Recall that on a given complete separable metric space $(S,\rho)$, the notion of weak convergence of sequences of Borel probability measures can be metrised by the bounded Lipschitz (BL) metric
\eq{
\label{Eq:BLMetric}
	d_S(\nu_1,\nu_2)=\sup_{F:S\to\R, \ \|F\|_{Lip}\le1}\left|\int_S F d(\nu_1-\nu_2)\right|, 
}
where
$$
	\|F\|_{Lip}=\sup_{x\in S}|F(x)|+\sup_{x,y\in S, \ x\neq y}\frac{|F(x)-F(y)|}{\rho(x,y)};
$$
see, e.g., \cite[Theorem 3.28]{Dudley2014}.

\begin{theorem}\label{Mtheo}

	Let $\Pi$ be a Gaussian Borel probability measure on $L^2(\OC)$ with RKHS $V_\Pi=H^r(\OC), \ r>d/2$. Assume that $f^\dagger\in H_c^\alpha(\OC), \ \alpha\ge r-d/2$, and let $P^M_{f^\dagger}$ be the law of $M_\varepsilon$ generated by (\ref{invprobl1}) with $f=f^\dagger$. Let $\bar f=\E^\Pi[f|M_\varepsilon]$ be the mean of the posterior distribution $\Pi(\cdot|M_\varepsilon)$ arising from observing \eqref{invprobl1}. Then, for all $\beta>2+d$ and any compact $K\subset\OC$, denoting $d_{(H^\beta_K)^*}$ the BL-metric for weak convergence on $(H^\beta_K(\OC))^* $, 
\eq{
\label{thesis1}
	d_{(H^\beta_K)^*}\left(\LC(\varepsilon^{-1}(f-\bar f)|M_\varepsilon),\mu\right)\to 0 
}
in $P^M_{f^\dagger}$-probability as $\varepsilon\to0$. Above $\LC(\varepsilon^{-1}(f-\bar f)|M_\varepsilon)$ is the centred and scaled posterior \eqref{centrScalPost}, and $\mu$ is the Gaussian distribution induced by X  in (\ref{X}).
\end{theorem}


	Similar results as Theorem \ref{thm:semiBvM2} and Theorem \ref{Mtheo} could be formulated for Example \ref{ex:psido}, exploiting the fact that the 'Fisher information operator' $A^*A$ has a well defined inverse $(A^*A)^{-1}:H^s(\Oo)\to H^{s-2t}(\Oo)$, for all $s\in\R$. In particular, since $\Oo$ was assumed to be a closed manifold, the weak convergence will be achieved in $H^{-\beta}(\Oo)$ for all $\beta>t+d$.

\begin{remark}[Applications to uncertainty quantification]\label{Remark:UQapplications}
	With similar reasoning as in Section \ref{Sec:UQ}, Theorem \ref{thm:semiBvM2} implies that for all $\psi\in H^\beta_c, \ \beta>2+d/2$, the credible intervals $C_\varepsilon$ in (\ref{Ceps}) centred at the plug-in Tikhonov regulariser $\langle \bar f,\psi\rangle_{L^2}$  have asymptotically correct frequentist coverage and optimal diameter.

	On the other hand, the full strength of  Theorem \ref{Mtheo} can be employed to show that the posterior distribution delivers valid uncertainty quantification also for the entire unknown $f$,  by considering credible sets in the weak topology where the limit is attained. The weak convergence in the dual space $(H^\beta_K)^*$ is indeed enough to deduce frequentist guarantees for a sufficiently rich class of credible sets (see the related discussion in Section 7.3.4 in \cite{Gine2016}). In particular, choosing posterior quantiles $\tilde R_\varepsilon=\tilde R(\alpha,M_\varepsilon)$ so that
$$
	\tilde C_\varepsilon=\{f\in L^2 : \|f-\bar f\|_{(H^\beta_K)^*}\le \tilde R_\varepsilon\},
	\quad \Pi(\tilde C_\varepsilon|M_\varepsilon)=1-\alpha, \ \alpha\in (0,1),
$$
we have for all $\beta>2+d$ 
$$
P^M_{f^\dagger}(f^\dagger\in \tilde C_\varepsilon)\to1-\alpha
$$
 as $\varepsilon\to0$, with asymptotically vanishing diameter $\tilde R_\varepsilon=O_{P^M_{f^\dagger}}(\varepsilon)$.

 	Finally, while the optimal rate $\varepsilon$ is obtained for the relatively weak norm of $(H^\beta_K)^*$, arguing as in Section 2 in \cite{Castillo2013} (see also Section 5.1 in \cite{Nickl2019}), we can intersect $\tilde C_\varepsilon$ with additional prior smoothness information (cf. Step I in Section \ref{sec:MtheoProof}) to show that the diameter of $\tilde C_\varepsilon$ decays at polynomial rate $\varepsilon^\gamma$, for any $\gamma<\alpha/(\alpha+2+d),$ also with respect to the stronger norm of interest $\|\cdot\|_{L^2(K')}$, for any compact $K'\subsetneq K$.
 	
\end{remark}

\begin{remark}[Smoothness requirement]\label{Rem:Smoothenss}
	Regarding the weak convergence to $\mu$ on $(H^\beta_K)^*$, the requirement that $\beta>2+d$ under which (\ref{thesis1}) is obtained is stronger than the necessary tightness condition $\beta>2+d/2$  of Lemma \ref{lemma:tightness}. While the proof of Theorem \ref{Mtheo} does imply the convergence of the finite-dimensional distributions of $\LC(\varepsilon^{-1}(f-\bar f)|M_\varepsilon)$ to those of $\mu$ in the full range $\beta>2+d/2$ (see Lemma \ref{lemma:fddconv}), the stronger condition $\beta>2+d$ is crucial in order to control the arising semiparametric bias term \textit{uniformly} in the collection $\{\psi\in H^\beta_K, \ \|\psi\|_{H^\beta}\le1\}$. This in turn implies that the $L^2$-diameter of $\tilde C_\varepsilon$ does not attain the minimax rate $\varepsilon^{\alpha/(\alpha+2+d/2)}$, which hence can potentially deliver polynomially sub-optimal result.

	To the best of our knowledge, examples of Gaussian priors that attain a nonparametric BvM theorem in the optimal function space are known in literature only in the SVD-based framework considered in  \cite{Castillo2013, Castillo2014,Ray2014thesis}, or in the 'nearly-diagonal' problem studied very recently by \cite{NicklRay2019}.
Applying our proof to a Gaussian prior defined via SVD would here recover the result of \cite{Ray2014thesis}.
However, the main interest of this paper is in the performance of standard Gaussian priors that are not defined on the SVD basis of the forward operator - such as the Mat\'ern priors considered in the examples - since this information is rarely available in inverse problems encountered in practice.
Our results show that for the inverse problem \eqref{invprobl1} standard Gaussian priors indeed yield optimal semiparametric inference for the maximal class of functionals, and provide a validation of the associated nonparametric credible sets.

\end{remark} 
	
%
%
%
%
%

\section{Proofs} 
\label{Sec:Proofs}  

%
%
%

\subsection{Proof of Theorem \ref{Thm:BvM1}.}
\label{Subsec:GeneralBvMProof}

	The proof of Theorem \ref{Thm:BvM1} follows ideas developed in \cite{Monard2019} for the special case of $A$ being the X-ray transform and $\psi\in C^\infty$. We will here outline the proof and comment on the main steps. We start by noting that the posterior concentrates on events that have high enough prior probability.  As a result, one can confine the analysis to an approximate posterior arising from restricting the prior over such sets. This observation allows to conveniently incorporate concentration properties of the prior into the analysis.

\begin{lemma}\label{lem:ae}

	Let $\Pi(\cdot \g M_\varepsilon)$ be the posterior distribution arising from observation $M_\varepsilon$ in \eqref{eq:observation} and prior $\Pi$ satisfying Condition \ref{cond:concentration} for a fixed $f^\dagger\in\Wt$ and some sequence $\delta_\varepsilon\to0$, such that ${\delta_\varepsilon}/{\varepsilon}\to\infty$. Then, for any Borel set $D_\varepsilon\subset\Wt$ for which 
\begin{align}
\label{eq:complementAs}
	\Pi(D_\varepsilon^c)
	\lesssim
		e^{-D(\delta_\varepsilon/\varepsilon)^2}, 
		\quad \text{for some $D>3$,}
\end{align}
and all $\varepsilon>0$ small enough, we have
\begin{align}
\label{Eq:Theo1Thesis2}
	\Pi(D_\varepsilon^c \g M_\varepsilon) \to 0 \quad \text{and}
	\quad 
	\|\Pi(\cdot \g M_\varepsilon)-\Pi^{D_\varepsilon}(\cdot \g M_\varepsilon)\|_{TV} \to 0
\end{align}
in $P_{f^\dagger}^M$-probability as $\varepsilon\to0$. Above $\Pi^{D_\varepsilon}(\cdot \g M_\varepsilon)$ is the posterior arising from the prior $\Pi(\cdot\cap D_\varepsilon)/\Pi(D_\varepsilon)$ restricted to $D_\varepsilon$ and renormalised.      

\end{lemma}

	The proof of Lemma \ref{lem:ae} (and Lemma \ref{lem:Convergence} below) can be adapted from the corresponding results in \cite{Monard2019}. They are included for completeness in the Supplement.

Next we need to find a suitable set $D_\varepsilon$. If $f\sim\Pi$, we have $\langle f,\pt\rangle_{V_\Pi}\sim \NC(0,\|\pt\|_{V_\Pi}^2)$ for all $\pt\in V_\Pi$, and the standard Gaussian tail bound guarantees for all $t\geq0$ that 
\begin{align*}
	\Pi\left(f : \frac{|\langle f,\pt\rangle_{V_\Pi}|}{\|\pt\|_{V_\Pi}}
	>\frac{t\delta_\varepsilon}{\varepsilon}\right)
	\leq
	e^{-\frac{t^2}{2}(\delta_\varepsilon/\varepsilon)^2}.
\end{align*} 
Hence we can choose     
\begin{align}
\label{simpleDeps}
	D_\varepsilon
	=
		\left\{f : \frac{|\langle f,\pt\rangle_{V_\Pi}|}{\|\pt\|_{V_\Pi}}
		\leq\frac{T\delta_\varepsilon}{\varepsilon}\right\},
	\quad T>\sqrt{6}.
\end{align}

We assume that the test function $\psi\in W_1$ fulfils $|\langle\psi,\varphi\rangle_{W_1}|\leq \|\varphi\|_\Wt$, for all $\varphi\in W_1$, in order to extend $L_\psi(\cdot)=\langle\psi,\cdot\rangle_{W_1}$ continuously to $\Wt$. If we assume furthermore that $\psi=-A^*A\pt$, with some $\pt\in V_\Pi$, we can proceed to study the moment generating function of $\varepsilon^{-1}(\langle f,\psi\rangle_{W_1}-\widehat{\Psi})$ under the posterior $\Pi^{D_\varepsilon}(\cdot\g M_\varepsilon)$, and conclude that it converges to the moment generating function of the limiting Gaussian law. 

\begin{lemma}\label{lem:Convergence}
	
	Under the conditions of Lemma \ref{lem:ae}, consider a test function $\psi\in W_1$ such that $|\langle\psi,\varphi\rangle_{W_1}|\lesssim \|\varphi\|_\Wt$, for all $\varphi\in W_1$,
and suppose that $\psi=-A^*A\pt$, for some $\pt \in V_\Pi$. Define the random variable
\begin{align*}
\widehat{\Psi}=\langle f^\dagger,\psi\rangle_{W_1}-\varepsilon\langle A\pt, \W\rangle_{W_2}. 
\end{align*}
Then, for all $\tau\in\R$ we have as $\varepsilon\to0$
\begin{align}
\label{eq:expected}
	\E^{\Pi^{D_\varepsilon}}\left(e^{\frac{\tau}{\varepsilon}
	(\langle f,\psi\rangle_{W_1}-\widehat{\Psi})}\Big{|} M_\varepsilon \right)
	=
		e^{\frac{\tau^2}{2}\|A\pt\|_{W_2}^2}\left(1+o_{P_{f^\dagger}^M}(1)\right).
\end{align}

\end{lemma}

	To conclude, we note that the exponential in the right hand side of \eqref{eq:expected} coincides with the moment generating function of $\NC(0,\|A\pt\|_{W_2}^2)$. Since the convergence of the Laplace transforms implies weak convergence (see, e.g., Proposition 1 in the supplement of \cite{Castillo2015}), we obtain from Lemma \ref{lem:Convergence} that the conclusion of Theorem \ref{Thm:BvM1} holds for the approximate posterior $\Pi^{D_\varepsilon}(\cdot|M_\varepsilon)$. Furthermore, convergence in total variation distance implies convergence in any metric for weak convergence and hence Theorem \ref{Thm:BvM1} follows from Lemma \ref{lem:ae}. 

%
%
%

\subsection{Proof of Proposition \ref{Ex:semiParBvM1}.}
\label{Sec:ProofExample}

	We now apply Theorem \ref{Thm:BvM1} to show the semiparametric result in the elliptic boundary value problem setting of Section \ref{sec:ellipticBVP}. As already noted before Proposition \ref{Ex:semiParBvM1}, any test function $\psi\in H^{r+4}_c=H^{r+4}_c(\OC)$ verifies the requirements of Theorem \ref{Thm:BvM1}. Hence, we only need to derive  Condition \ref{cond:concentration} for the chosen prior. In particular, for $\Pi$ a Gaussian prior on $L^2$ with RKHS $H^r$, $r>d/2$, and the true unknown $f^\dagger\in H^\alpha_c$, $\alpha\ge0$, we find suitable sequences $\delta_\varepsilon$ that satisfy the estimate \eqref{eq:speed} for the concentration function
\eq{\label{ConcFuncForL}
	\phi_{\Pi,f^\dagger}(\delta)
	=
		\inf_{g\in H^r, \ \|g-f^\dagger\|_{(H^2_0)^*}\le\delta}\frac{\|g\|^2_{H^r}}{2}
		-\log\Pi(f:\|f\|_{(H^2_0)^*}\le\delta),
	\quad \delta>0.
}
We proceed by calculating suitable upper bounds for the two terms. For the first term, we need to find approximations for the unknown $f^\dagger\in H^\alpha_c$ in the RKHS $H^r$ of $\Pi$, for which we can both control the approximation error and the norm in the latter space. We employ the approximations used in Section 4.3.3 of \cite{Nickl2007}. In particular, we fix a compact set $F$ such that $\supp(f^\dagger)\subsetneq F \subsetneq\OC$,  and a cut-off function $\zeta\in C^\infty_c(\R^d)$ such that $\zeta=1$ on $\supp(f^\dagger)$, $0\le\zeta\le1$  and $\supp(\zeta)\subseteq F.$ Noting that we can (isometrically) extend $f^\dagger$ to zero outside $F$ to form an element in $H^\alpha(\R^d)$, we then define
\eq{\label{EpsApprox}
f^\dagger_\varepsilon=(\zeta \FC^{-1}\1_{|\cdot|\le N_\varepsilon}\FC f^\dagger)|_{\OC},
}
for any sequence $N_\varepsilon\to\infty$ as $\varepsilon\to0$.

\begin{lemma}\label{lem:approximation}

	Let $f^\dagger\in H^\alpha_c(\OC)$ for some $\alpha>0$, and fix a compact set $F$ such that $\supp(f^\dagger)\subsetneq F\subsetneq \OC$. Then we have, for $f^\dagger_\varepsilon$ as in (\ref{EpsApprox}) and for any sequence $N_\varepsilon\to\infty$ as $\varepsilon\to0$, 

\begin{enumerate}

\item $f^\dagger_\varepsilon\in H_F^t(\OC)$ for all $t\ge0$ and
\eq{
\label{approxProp1}
	\|f^\dagger_\varepsilon\|^2_{H^t}
	\le 
		(1+N_\varepsilon^2)^{\max\{0,t-\alpha\}}\|f\|^2_{H^\alpha};
}

\item for all $0\le s<\alpha$
\eq{
\label{approxProp2}
	\|f^\dagger_\varepsilon-f^\dagger\|^2_{H^s}
	\le 
		(1+N^2_\varepsilon)^{s-\alpha}\|f^\dagger\|^2_{H^\alpha};
}
and for all $s\ge0$,
\eq{
\label{approxProp3}
	\|f^\dagger_\varepsilon-f^\dagger\|^2_{(H^s)^*}
	\le 
		(1+N^2_\varepsilon)^{-s-\alpha}\|f^\dagger\|^2_{H^\alpha}.
}

\end{enumerate}

\proof

	Let $t\ge0$ be fixed. Clearly $\supp(f^\dagger_\varepsilon)\subseteq\supp(\zeta)\subseteq F$, and we can compute directly 
\eqstar{
	 \| f^\dagger_\varepsilon\|^2_{H^t(\OC)} 
 	&\le
 		 \| \zeta \FC^{-1}\1_{|\cdot|\le N_\varepsilon}\FC f^\dagger\|^2_{H^t(\R^d)} \\
	&\lesssim
	 	\int_{\R^d}(1+|\xi|^2)^t(\1_{|\xi|\le N_\varepsilon}\FC f^\dagger(\xi))^2d\xi\\
	&=
	 	\int_{|\xi|\le N_\varepsilon}(1+|\xi|^2)^{t-\alpha}(1+|\xi|^2)^{\alpha}
		(\FC f^\dagger(\xi))^2d\xi\\
	&\le
	 	(1+N^2_\varepsilon)^{\max\{0,t-\alpha\}}\| f^\dagger\|^2_{H^\alpha(\OC)}.
}

	For $0\le s <\alpha$ we proceed similarly, observing that $f^\dagger=\zeta f^\dagger$ since $\zeta=1$ on $\supp(f^\dagger)$. Then
 
 \eqstar{
 	\| f^\dagger_\varepsilon -  f^\dagger\|^2_{H^s(\OC)}
	 &\le 
 		\| \zeta \FC^{-1}\1_{|\cdot|\le N_\varepsilon}
		\FC f^\dagger-\zeta f^\dagger\|^2_{H^s(\R^d)}\\
 	&\lesssim
 		\int_{\R^d}(1+|\xi|^2)^s(\1_{\{|\xi|\le N_\varepsilon\}}\FC f^\dagger(\xi)-
		\FC f^\dagger(\xi))^2d\xi\\
	& \le
	 	(1+N^2_\varepsilon)^{s-\alpha}\|f^\dagger\|^2_{H^\alpha(\OC)}.
 } 
 
	Finally, recalling that both $f^\dagger$ and $f^\dagger_\varepsilon$ are supported in $F\subsetneq\OC$, we get for all $s\ge0$
\eqstar{
	\|f^\dagger_\varepsilon-f^\dagger\|^2_{(H^s(\OC))^*}
	&=
		\sup_{u\in H^s(\OC), \|u\|_{H^s(\OC)}\le1}\Big{|}\int_F
		(f^\dagger_\varepsilon-f^\dagger)u dx\Big{|}\\
	&\le
		\sup_{U\in H^s(\R^d), \|U\|_{H^s(\R^d)}\le1}\Big{|}\int_F 
		(f^\dagger_\varepsilon-f^\dagger)U|_{\OC} dx\Big{|}\\
	&=
		\sup_{U\in H^s(\R^d), \|U\|_{H^s(\R^d)}\le1}\Big{|}\int_{F} 
		(\zeta\FC^{-1}\1_{|\cdot|\le N_\varepsilon}\FC f^\dagger- \zeta f^\dagger)U dx\Big{|}\\
	&\lesssim
		\|\FC^{-1}\1_{|\cdot|\le N_\varepsilon}\FC f^\dagger- f^\dagger\|_{H^{-s}(\R^d)}\\
	&\le
		(1+N_\varepsilon^2)^{-s-\alpha}\|f\|^2_{H^\alpha(\OC)},
}
where the last line follows arguing just as above for the case $0\le s<\alpha$.

\endproof

	We next derive an upper bound for the second term in \eqref{ConcFuncForL}. The proof adapts to the inverse problem \eqref{invprobl1} standard computations in the theory of small balls probabilities of Gaussian priors (e.g., \cite[Section 7.3]{Gine2016}).

\begin{lemma}\label{lem:SmallBallForL}

	Let $\Pi$ be a Gaussian Borel probability measure on $L^2(\OC)$ with RKHS $V_\Pi=H^r(\OC), \ r>d/2$. Then, as $\delta\to0$,
$$
	-\log \Pi( f : \|f\|_{(H_0^2(\OC))^*}\le \delta)\lesssim \delta^{-\frac{d}{r+2-d/2}}.
$$

\proof

	Since for any $f\in L^2$ we have $f=L(L^{-1}f)$, we can write
$$
	\Pi(f : \|f\|_{(H^2_0)^*}\leq\delta) = \Pi(f: \|L(L^{-1} f)\|_{(H_0^2)^*}\leq\delta).
$$

	Recalling that $L$ is self-adjoint when acting on $H^2_0$, we have for some $c>0$ that
\eqstar{
	\|L(L^{-1} f)\|_{(H^2_0)^*}
	=
		\sup_{v\in H^2_0,\ \|v\|_{H^2}\le 1}|\langle L(L^{-1}f),v\rangle_{L^2}| 
	\le
		c\| L^{-1}f\|_{L^2},
}
having used the boundedness of $L$. Thus,
\begin{align*}
\label{eq:lift}
	-\log \Pi(f : \|f\|_{(H^2_0)^*}\leq\delta) 
 	&\le 
		-\log \Pi\left(f : \|L^{-1} f\|_{L^2}\le\delta/c\right)\\
	&= 
		-\log\tilde\Pi \left(h : \|h\|_{L^2}\le\delta/c\right)
\end{align*}
where $h=L^{-1}f\sim \tilde\Pi$ for $f\sim\Pi$. From Exercise 2.6.5 in \cite{Gine2016} and the linearity of $L^{-1}$, we see that $\tilde\Pi$ is a Gaussian probability measure with RKHS $V_{\tilde\Pi}=L^{-1}(H^r)=H^{r+2}_0$, with unit ball $B_{V_{\tilde\Pi}}$ included in the unit ball $B^{r+2}$ of $H^{r+2}$. 
We thus get the following upper bound for the minimal number $N(B_{V_{\tilde\Pi}},\|\cdot\|_{L^2}, \delta)$ of $L^2$-balls or radius $\delta$ to cover $B_{V_{\tilde\Pi}}$:
\begin{align*}
	N(B_{V_{\tilde\Pi}},\|\cdot\|_{L^2}, \delta) 
	\le
	 	N(B^{r+2},\|\cdot\|_{L^2}, \delta).
\end{align*}
Theorem 4.3.36 in \cite{Gine2016} now implies that
$$  
	\log N(B^{r+2},\|\cdot\|_{L^2}, \delta) 
	\lesssim
		\delta^{-\frac{d}{r+2}},
$$
and by applying the small ball estimates in Theorem 1.2 of \cite{Li1999}, we obtain that as $\delta\to0$
$$
-\log\tilde\Pi\left(h : \|h\|_{L^2}\le\delta/c\right)\lesssim \delta ^{-\frac{d}{r+2-d/2}},
$$
concluding the proof.

\endproof

\end{lemma}

\end{lemma}

	Thus, applying Lemma \ref{lem:approximation}, the first term in the estimate (\ref{ConcFunAsymptForL}) follows by choosing, for any fixed $\delta\ge0$, $N_\varepsilon$ in (\ref{EpsApprox}) in such a way that
$$
(1+N_\varepsilon^2)^{-2-\alpha}\le \delta^2,
$$
so that, in view of  (\ref{approxProp2}) and (\ref{approxProp1}) respectively,
$
\|f^\dagger_\varepsilon-f^\dagger\|_{(H^2_0)^*}\lesssim \delta
$
and
$
\|f^\dagger_\varepsilon\|_{H^r}^2 
\lesssim 
 	\delta^{-\frac{2\max\{0,r-\alpha\}}{2+\alpha}}.
$
It can then be readily checked from (\ref{ConcFunAsymptForL}) that the sequence $\delta_\varepsilon$ in (\ref{DeltaEpsForL}) satisfies the required inequality 
$
	\phi_{\Pi,f^\dagger}(\delta_\varepsilon)
	\lesssim 
		(\delta_\varepsilon/\varepsilon)^2,
$
concluding the proof of Proposition \ref{Ex:semiParBvM1}.

%
%
%

\subsection{Proofs of Theorem \ref{thm:semiBvM2} and Theorem \ref{Mtheo}.}
\label{sec:MtheoProof}

	The key steps of the proof consist in a refinement of the strategy developed to prove Theorem \ref{Thm:BvM1} and Proposition \ref{Ex:semiParBvM1}. Following \cite{Castillo2013,Castillo2014,Nickl2017}, we first aim at obtaining the Laplace transform convergence \eqref{eq:expected} uniformly with respect to the test functions $\psi\in H^\beta_c, \ \beta>2+d/2$ (cf. Steps I-II). We subsequently exploit the result to show Theorem \ref{thm:semiBvM2}, and to derive the convergence of the finite dimensional distributions of the centred and scaled posterior $\LC(\varepsilon^{-1}(f-\bar f)|M_\varepsilon)$ to those of the limiting Gaussian measure $\mu$ (Step III). Finally, combining this observation with a suitable bound on the covariance of the process $X_\varepsilon$ in (\ref{postsp}), we show for each $\beta>2+d$ that the distance between $\LC(\varepsilon^{-1}(f-\bar f)|M_\varepsilon)$ and $\mu$, measured in the BL metric on $(H^\beta_K)^*$, vanishes with $P_{f^\dagger}^M$-probability converging to one (Step IV-V).

\paragraph{\bf Step I: Construction of the approximating sets.}
\label{Par:ApproxSets} 

	Let $\Pi$ be a centred Gaussian prior on $L^2=L^2(\OC)$ with RKHS $H^r$, $r>d/2$, and let $f^\dagger\in H_c^\alpha$ be fixed. Recall that we assume the prior to undersmooth $f^\dagger$, namely that $\alpha\geq r-d/2$. Then Remark \ref{Rem:ConcFunAsymptForL} implies that Condition \ref{cond:concentration} is satisfied by  taking
\eq{
\label{Eq:DeltaEps}
\delta_\varepsilon\simeq \varepsilon^{\frac{2+r-d/2}{2+r}}.
}

	In the first step we need to construct appropriate approximating sets $D_\varepsilon$, by adapting the events introduced in (\ref{simpleDeps}) for the proof of Theorem \ref{Thm:BvM1}. First, to extend the semiparametric result in Proposition \ref{Ex:semiParBvM1} to the range $2+d/2<\beta<r+4$, we replace the element $\tilde\psi=-LL\psi$ (here not in the RKHS $V_\Pi=H^r$) with a suitable approximation. To deal with the possibly diverging norm of such approximations, we will then impose additional constraints to control the size of $f\in D_\varepsilon$. Finally, to achieve the required uniformity in the Laplace transform convergence \eqref{eq:expected}, we will further intersect the resulting events across all test functions, in such a way as to maintain the exponential decay \eqref{eq:complementAs} for $\Pi(D_\varepsilon^c)$.

	To proceed, let $\beta>2+d/2$, let $K\subset\OC$ be compact and fix a compact set $F$ such that $K\subsetneq F\subsetneq \OC$. Then, for each $\psi$ in a ball 
\begin{equation}
\label{Eq:Balls}
B^\beta_K(z):=\{v\in H_K^\beta, \|v\|_{H^\beta}\le z\}
\end{equation} 
of fixed radius $z>0$, consider the approximation of $L\psi$ given by Lemma \ref{lem:approximation}, of the form
\eq{
\label{tildepsieps}
	\tilde\psi_\varepsilon
	=
		(\zeta \FC^{-1}\1_{|\cdot|\le N_\varepsilon}\FC[ L\psi])|_{\OC},
		\quad
		\textnormal{with}\ \ N_\varepsilon\simeq \varepsilon^{-\frac{1}{2+r}}.
}
By point 1. in Lemma \ref{lem:approximation}, we can uniformly control the Sobolev norms of the resulting collection of approximations. Indeed, by the continuity of $L$, for all $\psi\in B^\beta_K(z)$ we have $\|L\psi\|_{H^{\beta-2}}\le z'$ for some constant $z'>0$, so that in view (\ref{approxProp1}), for all $t\ge0$,
\eq{\label{CollectionTildePsiEps}
	\{\tilde\psi_\varepsilon, \ \psi\in B^\beta_K(z)\}
	\subseteq B_F^t(b^t_\varepsilon),
	 \quad 
	 b^t_\varepsilon
	 	:=
			\sup_{ \psi\in B^\beta_K(z)}\|\tilde\psi_\varepsilon\|_{H^t}					\le
			z'(1+N^2_\varepsilon)^{\max\{0,t-\beta+2\}/2}.
}
Then, for all $\psi\in B^\beta_K(z)$, it follows in particular $\tilde\psi_\varepsilon\in B_F^{r+2}( b^{r+2}_\varepsilon)$, from which we deduce that $L\tilde\psi_\varepsilon\in V_\Pi= H^r$. Thus, if $f\sim\Pi$, then $\langle f,L\tilde\psi_\varepsilon \rangle_{V_\Pi}\sim\NC(0,\|L\tilde\psi_\varepsilon\|^2_{V_\Pi})$, with variance uniformly bounded, in view of the isomorphism property of $L$, by
\eq{
\label{sigmaepsilon}
	\sigma^2_\varepsilon
	:=
		\sup_{\psi\in B^\beta_K(z)}\E^\Pi|\langle f,
		L\tilde\psi_\varepsilon \rangle_{V_\Pi}|^2
	\simeq 
		(b^{r+2}_\varepsilon)^2.
}

	Define, for each $\varepsilon>0$, and $D>0$ to be chosen below, the approximating set
\eq{
\label{setGeps}
G_\varepsilon
=
	\left\{f:\sup_{\psi\in B^\beta_K(z)}|\langle f,L\tilde\psi_\varepsilon \rangle_{V_\Pi}|\le 
	D\sigma_\varepsilon\delta_\varepsilon/\varepsilon\right\}.
}
Here $G_\varepsilon$ serves as the counterpart of the events (\ref{simpleDeps}), with the constraint holding simultaneously for all $\psi$. 

	We derive the exponential decay \eqref{eq:complementAs} for $\Pi(G_\varepsilon^c)$. First, denoting $\E^\Pi$ the expectation under the prior, we have by the Borell-Sudakov-Tirelson inequality \cite[Theorem 2.5.8]{Gine2016} that for all $\tilde D>0$
\eq{
\label{Eq:BorTsirSudIneq}
	\Pi\left(f:\sup_{\psi\in B^\beta_K(z)}|\langle f,L\tilde\psi_\varepsilon \rangle_{V_\Pi}|
	> \E^\Pi\sup_{\psi\in B^\beta_K(z)}|\langle f,L\tilde\psi_\varepsilon \rangle_{V_\Pi}| 
	+ \tilde D\sigma_\varepsilon\delta_\varepsilon/\varepsilon\right)
	\le 
		e^{-\frac{\tilde D^2}{2}(\delta_\varepsilon/\varepsilon)^2}.
}
Thus, the condition \eqref{eq:complementAs} will follow if we show that 
$
	\E^\Pi\sup_{\psi\in B^\beta_K(z)}|\langle f,L\tilde\psi_\varepsilon \rangle_{V_\Pi}|
	\lesssim
		\sigma_\varepsilon\delta_\varepsilon/\varepsilon.
$ 
Indeed, in view of (\ref{CollectionTildePsiEps}), denoting $B^s(z)$ a ball in $H^s$ of radius $z$, for general $s\ge0$ and $z>0$,
$$
\E^\Pi\sup_{\psi\in B^\beta_K(z)}|\langle f,L\tilde\psi_\varepsilon \rangle_{V_\Pi}|
\le
	\E^\Pi\sup_{v\in B^{t+2}(z'b^{t+2}_\varepsilon)}|\langle f,Lv\rangle_{V_\Pi}|
\lesssim
	\E^\Pi\sup_{w\in B^t(z''b^{t+2}_\varepsilon)}|\langle f,w\rangle_{V_\Pi}|
$$
and Dudley's bound for the expectation of suprema of Gaussian processes \cite[Theorem  2.3.8]{Gine2016} yields, for $\sigma_\varepsilon$ the constant in (\ref{sigmaepsilon}),
\eqstar{
	\E^\Pi\sup_{w\in B^t(z''b^{t+2}_\varepsilon)}|\langle f,w\rangle_{V_\Pi}| 
	&\lesssim
		\int_0^{\sigma_\varepsilon}
		\sqrt{\log N(\eta,B^t(z''b^{t+2}_\varepsilon),\|\cdot\|_{V_\Pi})}d\eta\\
	&=
		\int_0^{\sigma_\varepsilon}
		\sqrt{\log N\left(\frac{\eta}{z''b^{t+2}_\varepsilon},B^t(1),\|\cdot\|_{V_\Pi}\right)}d\eta.
}
Fixing $t>\max\{r+d/2,\beta-4\}$, recalling $V_\Pi=H^r$ and using the known metric entropy estimates for Sobolev balls (see, e.g., \cite{Triebel1978}), we then obtain
\eqstar{
	\E^\Pi\sup_{w\in B^t(z''b^{t+2}_\varepsilon)}|\langle f,w\rangle_{V_\Pi}| 
	&\lesssim 
		\int_0^{\sigma_\varepsilon}
		\left(b^{t+2}_\varepsilon/\eta\right)^{\frac{d}{2(t-r)}}d\eta
	\lesssim
		(b^{t+2}_\varepsilon)^{\frac{d}{2(t-r)}} \sigma_\varepsilon^{\frac{2t-2r-d}{2(t-r)}}.
}
Using (\ref{sigmaepsilon}),  it follows that
\eqstar{
	\frac{1}{\sigma_\varepsilon}
	\E^\Pi\sup_{w\in B^t(z''b^{t+2}_\varepsilon)}|\langle f,w\rangle_{V_\Pi}|
	&\lesssim
		(b^{t+2}_\varepsilon)^{\frac{d}{2(t-r)}} (b^r_\varepsilon)^{-\frac{d}{2(t-r)}}\\
	&= 
		(1+N^2_\varepsilon)^{\frac{d(t-\beta+4)}{4(t-r)}} 
		(1+N^2_\varepsilon)^{-\frac{d\max\{0,r-\beta+4\}}{4(t-r)}}\\
	&\le
		(1+N^2_\varepsilon)^{d/4}.
}
Recalling that $\delta_\varepsilon\simeq\varepsilon^{\frac{2+r-d/2}{2+r}}$ the choice $N_\varepsilon\simeq \varepsilon^{-\frac{1}{2+r}}$, finally yields
\eqstar{
	\frac{1}{\sigma_\varepsilon}
	\E^\Pi\sup_{\psi\in B^\beta_K(z)}|\langle f,L\tilde\psi_\varepsilon \rangle_{V_\Pi}| 
	\lesssim 
		 \varepsilon^{-\frac{d}{2(2+r)}}
	\simeq
		 \delta_\varepsilon/\varepsilon.
}
Taking $\tilde D>\sqrt{6}$ in \eqref{Eq:BorTsirSudIneq}, and sufficiently large $D>\tilde D$ in the definition (\ref{setGeps}) of $G_\varepsilon$, yields the exponential decay \eqref{eq:complementAs} for $\Pi(G_\varepsilon^c)$.

	Next, we proceed by suitably controlling the size of the elements in the approximating sets. To do so, let, for $\Phi$ is the standard normal cumulative distribution function,
$$
	Q_\varepsilon
	=
		-2\Phi^{-1}\left(e^{-\frac{\tilde D^2}{2}(\delta_\varepsilon/\varepsilon)^2}\right)
	\simeq
		\delta_\varepsilon/\varepsilon.
$$
For $\rho>0$ to be chose below and arbitrary $\kappa>0$, Consider the event
\eq{
\label{Feps}
	F_\varepsilon
	=
		\{f=f_1+f_2: \|f_1\|_{L^2}\le \rho\varepsilon^{\frac{r-d/2}{2+r}},
		\ \|f_2\|_{V_\Pi}\le Q_\varepsilon+\kappa\}
}
in which we constraint the prior draws $f\sim\Pi$ to belong to (a slight enlargement of) a ball of the RKHS $V_\Pi$ of growing radius. By the isoperimetric inequality  for Gaussian processes \cite[Theorem 2.6.12]{Gine2016} we can lower bound the prior probability of $F_\varepsilon$ by
\eq{
\label{eq11}
	\Pi(F_\varepsilon) 
	&\ge
		\Phi(\Phi^{-1}[\Pi(f:\|f\|_{L^2}\le \rho\varepsilon^{\frac{r-d/2}{2+r}})]+Q_\varepsilon),
}
Applying again the small ball estimate for $\Pi$ in Theorem 1.2 in \cite{Li1999} as in the proof of Lemma \ref{lem:SmallBallForL}, we see that for some $b>0$
$$
-\log\Pi(f:\|f\|_{L^2}\le\rho\varepsilon^{\frac{r-d/2}{2+r}}) \le b\rho^{-\frac{d}{r-d/2}}\varepsilon^{-\frac{d}{2+r}}
$$
and recalling that $\delta_\varepsilon/\varepsilon\simeq\varepsilon^{-\frac{d/2}{2+r}}$, we can choose $\rho>0$ so that
\eqstar{
	-\log\Pi(f:\|f\|_{L^2}
	\le
		\rho\varepsilon^{\frac{r-d/2}{2+r}})
	\le
		 (\delta_\varepsilon/\varepsilon)^2.
}
Combining the above with (\ref{eq11}) yields
\eqstar{
	\Pi(F_\varepsilon) 
	&\ge
		 \Phi(\Phi^{-1}( e^{-(\delta_\varepsilon/\varepsilon)^2})+Q_\varepsilon)
	 \ge
	 	 \Phi(\Phi^{-1}(e^{-\frac{\tilde D^2}{2}(\delta_\varepsilon/\varepsilon)^2})+Q_\varepsilon)
	=
		 \Phi(Q_\varepsilon/2),
}
and finally 
$
	\Pi(F_\varepsilon^c)
	\le 
		e^{-\frac{\tilde D^2}{2}(\delta_\varepsilon/\varepsilon)^2}.
$

	We conclude by taking
\eq{\label{newDeps}
D_\varepsilon=G_\varepsilon\cap F_\varepsilon,
}
for which the bounds on $\Pi(G^c_\varepsilon)$ and $\Pi(F^c_\varepsilon)$ imply
$
\Pi(D_\varepsilon^c)
\le 2e^{-\frac{\tilde D^2}{2}(\delta_\varepsilon/\varepsilon)^2},\ 
$
$\tilde D^2/2>3$.

\paragraph{\bf Step II: Laplace transform expansion.}
\label{Step:LaplTransExp}
 
	We proceed deriving an asymptotic expression, analogous to the one obtained in Lemma \ref{lem:Convergence}, for  the Laplace transform of the linear functionals $\langle f,\psi\rangle_{L^2}$. In view of the simultaneous constraint imposed in \eqref{setGeps}, the result holds uniformly with respect to test functions $\psi$.

\begin{lemma}\label{Laplexpansion}

	Let $\Pi$ be a Gaussian Borel probability measure on $L^2(\OC)$ with RKHS $V_\Pi=H^r(\OC), \ r>d/2$, and assume that $f^\dagger\in H_c^\alpha(\OC), \ \alpha\ge r-d/2$.
For all $\beta>2+d/2$, and any $\psi\in B_K^\beta(z),\ z>0,$ (defined as in \eqref{Eq:Balls}) let $\tilde\psi_\varepsilon$ be the approximation in (\ref{tildepsieps}), and define
\eq{
\label{Eq:HatPsiPsi}
	\hat\Psi(\psi)
	=
		\langle f^\dagger,\psi\rangle_{L^2}
		+\varepsilon\langle \tilde\psi_\varepsilon,\W\rangle_{L^2}.
}
Then, for all fixed $\tau\in\R$
\eq{\label{finalconv2}
	\E^{\Pi^{D_\varepsilon}}\left[e^{\frac{\tau^2}{\varepsilon}
	[\langle f,\psi\rangle_{L^2}-\hat\Psi(\psi)]}\Big{|}M_\varepsilon\right]
	=
		e^{R_\varepsilon}e^{\frac{\tau^2}{2}\|L\psi\|^2_{L^2}}
		\frac{\Pi(D_{\varepsilon,\tau}|M_\varepsilon)}{\Pi(D_{\varepsilon}|M_\varepsilon)},
}   
where $D_{\varepsilon,\tau}=\{f-\tau\varepsilon L\tilde\psi_\varepsilon,\ f\in D_\varepsilon\}$ and $R_\varepsilon\to 0$  uniformly in $B^\beta_K(z)$ for any $z>0$ as $\varepsilon\to0$.

\proof
    
	We have
$$
	\E^{\Pi^{D_\varepsilon}}\left[e^{\frac{\tau^2}{\varepsilon}
	[\langle f,\psi\rangle_{L^2}-\hat\Psi(\psi)]}\Big{|}M_\varepsilon\right]
	=
		 e^{-\tau\langle\tilde\psi_\varepsilon,\W\rangle_{L^2}}
		 \E^{\Pi^{D_\varepsilon}}
		 \left[e^{\frac{\tau}{\varepsilon} \langle f-f^\dagger,\psi\rangle_{L^2}}
		 \Big{|}M_\varepsilon\right],
$$
and letting $f_{\tau}=f-\tau \varepsilon L\tilde\psi_\varepsilon$, the expectation in the right hand side becomes (cf. \eqref{eq:BayesFormula})
$$
	\E^{\Pi^{D_\varepsilon}}
	 \left[e^{\frac{\tau}{\varepsilon} \langle f-f^\dagger,\psi\rangle_{L^2}}
	 \Big{|}M_\varepsilon\right]
	 = 
		 \frac{\int_{D_\varepsilon}
		 e^{\frac{\tau}{\varepsilon}\langle f-f^\dagger,\psi\rangle_{L^2}}
		 e^{\ell(f)-\ell(f_\tau)}e^{\ell(f_\tau)}d\Pi(f)}{\int_{D_\varepsilon}e^{\ell(f)}d\Pi(f)}.
$$
From the expression of the log-likelihood \eqref{Likelihood} we readily obtain 
\eqstar{
	\ell(f)-\ell(f_\tau)
	&=
		 \frac{\tau^2}{2}\|\tilde\psi_\varepsilon\|^2_{L^2}
		 +\frac{\tau}{\varepsilon}
		 \langle L^{-1}(f-f^\dagger),\tilde\psi_\varepsilon\rangle_{L^2}
		 +\tau\langle\tilde\psi_\varepsilon, \W \rangle_{L^2},
}
which substituted into the previous expression yields, using the self-adjointness of $L^{-1}$,
\eq{
\label{lapl3}
	\E^{\Pi^{D_\varepsilon}}&\left[e^{\frac{\tau^2}{\varepsilon}
	[\langle f,\psi\rangle_{L^2}-\hat\Psi(\psi)]}\Big{|}M_\varepsilon\right]\\
	&=
		e^{\frac{\tau^2}{2}\|\tilde\psi_\varepsilon\|^2_{L^2}}
		e^{-\frac{\tau}{\varepsilon} 
		\langle L^{-1}f^\dagger,L\psi-\tilde\psi_\varepsilon\rangle_{L^2}}
		\frac{\int_{D_\varepsilon}e^{\frac{\tau}{\varepsilon}\langle L^{-1}f,
		L\psi-\tilde\psi_\varepsilon\rangle_{L^2}}
		e^{\ell(f_\tau)}d\Pi(f)}{\int_{D_\varepsilon}e^{\ell(f)}d\Pi(f)}.
}

	In view of \eqref{approxProp2}, we have that $\|\tilde\psi_\varepsilon-L\psi\|_{L^2}\to0$ as $\varepsilon\to0$ uniformly in $B^\beta_K(z)$ for all $z>0$, and hence
\eq{
\label{canc1}
	e^{\frac{\tau^2}{2}\|\tilde\psi_\varepsilon\|^2_{L^2}}
	=
		(1+o(1))e^{\frac{\tau^2}{2}\|L\psi\|^2_{L^2}}.
}
Next we prove that
\eq{
\label{canc2}
	e^{-\frac{\tau}{\varepsilon} 
	\langle L^{-1}f^\dagger,L\psi-\tilde\psi_\varepsilon\rangle_{L^2}}
	=
		1+o(1)
	\quad 
	\textnormal{uniformly in $B^\beta_K(z).$}
}
To do so, notice 
\eqstar{
	\sup_{\psi\in B^\beta_K(z)}
	&\bigg{|}-\frac{1}{\varepsilon}\langle L^{-1}f^\dagger,L\psi-
	\tilde\psi_\varepsilon\rangle_{L^2}\bigg{|}\\
	&=
		\frac{1}{\varepsilon}\| L^{-1}f^\dagger\|_{H^{\alpha+2}}
		\sup_{\psi\in B^\beta_K(z)}\left| 	
		\left\langle \frac{L^{-1}f^\dagger}{\| L^{-1}f^\dagger\|_{H^{\alpha+2}}},
		L\psi-\tilde\psi_\varepsilon\right\rangle_{L^2}\right|\\
	&\lesssim 
		\frac{1}{\varepsilon}
		\sup_{\psi\in B^\beta_K(z)}
		\|L\psi-\tilde\psi_\varepsilon\|_{(H^{\alpha+2}_0)^*}\\
	&\lesssim 
		\varepsilon^{-1}(1+N_\varepsilon^2)^\frac{-\alpha-\beta}{2}
}
where the last line follows by (\ref{approxProp3}). Recalling that $N_\varepsilon\simeq\varepsilon^{-\frac{1}{2+r}}$, $\alpha\ge r-d/2$ and $\beta>2+d/2,$
$$
	\sup_{\psi\in B^\beta_K(z)}\bigg{|}-\frac{\tau}{\varepsilon}
	\langle L^{-1} f^\dagger,L\psi-\tilde\psi_\varepsilon\rangle_{L^2}\bigg{|}
	\lesssim
		\frac{1}{\varepsilon}(1+N_\varepsilon^2)^{\frac{-\alpha-\beta}{2}}
	\simeq
		\varepsilon^\frac{\alpha+\beta-2-r}{2+r}
	\to
		0.
$$

	The following step consists in showing that uniformly in $B^\beta_K(z)$
\eq{
	\label{canc3}
	\int_{D_\varepsilon}e^{\frac{\tau}{\varepsilon}\langle L^{-1}f,
	L\psi-\tilde\psi_\varepsilon\rangle_{L^2}}e^{\ell(f_\tau)}d\Pi(f)
	=
		(1+o(1))\int_{D_\varepsilon}e^{\ell(f_\tau)}d\Pi(f).
}
The result will follow from the dominated convergence theorem upon showing that
\eqstar{
\sup_{f\in D_\varepsilon}\sup_{\psi\in B^\beta_K(z)}&\left|\frac{\tau}{\varepsilon}\langle L^{-1} f,L\psi-\tilde\psi_\varepsilon\rangle_{L^2}\right| \to0.
}
 Recalling the definition (\ref{newDeps}) of $D_\varepsilon$, we bound the left hand side by
\eqstar{
	\sup_{f\in F_\varepsilon}
	&\sup_{\psi\in B^\beta_K(z)}\left|\frac{\tau}{\varepsilon}
	\langle L^{-1}f,L\psi-\tilde\psi_\varepsilon\rangle_{L^2}\right| \\
	&\lesssim  
		\frac{1}{\varepsilon}\sup_{\|f_1\|_{L^2}\le \rho\varepsilon^{\frac{r-d/2}{2+r}}}
		\sup_{\psi\in B^\beta_K(z)}\left|\langle L^{-1}f_1,L\psi-							\tilde\psi_\varepsilon\rangle_{L^2}\right|\\
	&\quad
		+ 
		\frac{1}{\varepsilon}\sup_{\|f_2\|_{V_\Pi}\le Q_\varepsilon+\kappa}
		\sup_{\psi\in B^\beta_K(z)}\left|\langle L^{-1}f_2,
		L\psi-\tilde\psi_\varepsilon\rangle_{L^2}\right|.
}
Accordingly, it is enough to show the joint convergence of the two terms above, which can be done similarly as in the derivation of (\ref{canc2}). In particular
\eqstar{
 	\frac{1}{\varepsilon}
	\sup_{\|f_1\|_{L^2}\le \rho\varepsilon^{\frac{r-d/2}{2+r}}}
	\sup_{\psi\in B^\beta_K(z)}\left|\langle L^{-1}f,
	L\psi-\tilde\psi_\varepsilon\rangle_{L^2}\right|
	& \lesssim
		\varepsilon^{\frac{r-d/2}{2+r}-1}(1+N^2_\varepsilon)^{-\frac{\beta}{2}}
	\to
		0.
} 
On the other hand, recalling that $V_\Pi=H^r$, 
\eqstar{
	 \frac{1}{\varepsilon}\sup_{\|f_2\|_{V_\Pi}\le Q_\varepsilon+\kappa}
	 \sup_{\psi\in B^\beta_K(z)}\left|\langle L^{-1}f_2,
	 L\psi-\tilde\psi_\varepsilon\rangle_{L^2}\right|
	& \lesssim
		\varepsilon^{-1}(Q_\varepsilon+\kappa)(1+N_\varepsilon^2)^{-\frac{r+\beta}{2}}
	\to
		0,
}
since $\delta_\varepsilon\simeq\varepsilon^{\frac{2+r-d/2}{2+r}}$ and $
Q_\varepsilon=-2\Phi^{-1}\left(c'e^{-\frac{\tilde D^2}{2}(\delta_\varepsilon/\varepsilon)^2}\right)\simeq \delta_\varepsilon/\varepsilon$. Replacing (\ref{canc1}), (\ref{canc2}) and (\ref{canc3}) into (\ref{lapl3}) we obtain, uniformly in $B^\beta_K(z)$,
\eq{\label{lapl4}
	\E^{\Pi^{D_\varepsilon}}\left[e^{\frac{\tau}{\varepsilon}[\langle f,
	\psi\rangle_{L^2}-\hat\Psi(\psi)]}\big{|}M_\varepsilon\right] 
	=
		 (1+o(1)) 
		 e^{\frac{\tau^2}{2}\|L\psi\|^2_{L^2}}\frac{\int_{D_\varepsilon}e^{\ell(f_\tau)}d\Pi(f)}
		 {\int_{D_\varepsilon}e^{\ell(f)}d\Pi(f)}.
}

	We conclude by further simplifying the ratio in the right hand side of (\ref{lapl4}) in the same way as in the conclusion of the proof of Proposition 3.2 in \cite{Monard2019}. Let $\Pi_{\tau}$ be the law of the shifted parameter $f_{\tau}=f-\tau \varepsilon L\tilde\psi_\varepsilon$,  
and $D_{\varepsilon,\tau}=\{f-\tau \varepsilon L\tilde\psi_\varepsilon,\ f\in D_\varepsilon\}$. Then, the Cameron-Martin theorem (e.g., Theorem 2.6.13 in \cite{Gine2016}) yields
\eqstar{
	\frac{\int_{D_\varepsilon}e^{\ell(f_\tau)}d\Pi(f)}{\int_{D_\varepsilon}e^{\ell(f)}d\Pi(f)}
	&= 
		e^{-\frac{(\tau\varepsilon)^2}{2}\|L\tilde\psi_\varepsilon\|^2_{V_\Pi}}
		\frac{\int_{D_{\varepsilon,\tau}}e^{\ell(g)}										e^{-\tau\varepsilon\langle L\tilde\psi_\varepsilon,g\rangle_{V_\Pi}}d\Pi(g)}				{\int_{D_\varepsilon}e^{\ell(g)}d\Pi(g)}.
}
First notice that by (\ref{approxProp1}),
\eq{
\label{cancx}
	\sup_{\psi\in B^\beta_K(z)} \varepsilon^2\|L\tilde\psi_\varepsilon\|^2_{V_\Pi} 
	\lesssim 
		\sup_{\psi\in B^\beta_K(z)} \varepsilon^2\|\tilde\psi_\varepsilon\|^2_{H^{r+2}} 
	\lesssim 
		\varepsilon^{2-\frac{2\max\{0,r-\beta+4\}}{2+r}}\to0.
}
Next, recalling the definitions (\ref{setGeps}) and (\ref{newDeps}) of $G_\varepsilon$ and $D_\varepsilon$ respectively, we have
\eqstar{
	\sup_{g\in D_{\varepsilon,\tau}}\sup_{\psi\in B^\beta_K(z)} 
	|\varepsilon\langle L\tilde\psi_\varepsilon,g\rangle_{V_\Pi}|
	&=
		 \varepsilon \sup_{f\in D_{\varepsilon}} \sup_{\psi\in B^\beta_K(z)}
		|(L\tilde\psi_\varepsilon,f-\tau \varepsilon L\tilde\psi_\varepsilon\rangle_{V_\Pi}|\\
	&\le
		 \varepsilon\left(\sup_{f\in G_{\varepsilon}}\sup_{\psi\in B^\beta_K(z)}|				\langle L\tilde\psi_\varepsilon,f\rangle_{V_\Pi}|+|\tau|\varepsilon\sup_{\psi\in 			B^\beta_K(z)}\|	L\tilde\psi_\varepsilon\|_{V_\Pi}^2\right)\\
&\lesssim
	 \sigma_\varepsilon\delta_\varepsilon+o(1)
}
and, since $\beta>2+d/2$ and $r> d/2$, then
$
	\sigma_\varepsilon\delta_\varepsilon
	\lesssim
	\varepsilon^{\frac{\min\{r+2-d/2,\beta-2-d/2\}}{2+r}}
	\to
		0.
$
Thus, we conclude that uniformly in $B^\beta_K(z)$ as $\varepsilon\to0$,
\eqstar{
	\frac{\int_{D_\varepsilon}e^{\ell(f_\tau)}d\Pi(f)}{\int_{D_\varepsilon}e^{\ell(f)}d\Pi(f)}
	=
		(1+o(1))\frac{\int_{D_{\varepsilon,\tau}}e^{\ell(g)}d\Pi(g)}
		{\int_{D_\varepsilon}e^{\ell(g)}d\Pi(g)}
	=(1+o(1))
		\frac{\Pi(D_{\varepsilon,\tau}|M_\varepsilon)}{\Pi(D_{\varepsilon}|M_\varepsilon)}.
}

\endproof

\end{lemma}

\paragraph{\bf Step III: Convergence of the finite dimensional distributions.}

	We now exploit the previous lemma to show the convergence of the finite dimensional distributions of the centred and scaled posterior to those of the Gaussian measure $\mu$  induced by the process $X$ in (\ref{X}). The result is obtained by showing that the ratio in the right hand side of \eqref{finalconv2} converges to one, yielding, as in Lemma \ref{lem:Convergence}, the desired asymptotic expression for the Laplace transform. This will in turn conclude the proof of Theorem \ref{thm:semiBvM2}.

	To proceed, consider the centring Gaussian process obtained by collecting the random variables introduced in \eqref{Eq:HatPsiPsi},
\eq{
\label{hatpsieps}
	\hat \Psi_\varepsilon
	=
		(\hat\Psi_\varepsilon(\psi) : \psi\in H^\beta_c),
		 \quad \hat\Psi_\varepsilon(\psi)
	=
		\langle f^\dagger,\psi\rangle_{L^2}+\varepsilon\langle \tilde\psi_\varepsilon,
		\W\rangle_{L^2}, \quad \varepsilon>0.
}
In view of Lemma \ref{lemma:tightness}, and by the continuity of the linear map 
$
\psi\in H^\beta_c\mapsto \tilde\psi_\varepsilon\in H^t, \ t\ge0,
$
$\hat \Psi_\varepsilon$ defines a Borel measurable map on $(H^\beta_c)^*$, for $\beta>2+d/2$. Then, denote by 
\eq{
\label{almostPost}
	\LC(\varepsilon^{-1}(f-\hat\Psi_\varepsilon)|M_\varepsilon)=\LC(\hat X_\varepsilon), \quad 	f\sim\Pi
}
the tight condition law on $(H^\beta_c)^*$ of
\eq{
\label{hatXeps}
	\hat X_\varepsilon
	=
		(\hat X_\varepsilon(\psi) : \psi\in H^\beta_c), 
	\quad
			\hat 	X_\varepsilon(\psi)
	=
			\varepsilon^{-1}[\langle f,\psi\rangle_{L^2}-\hat\Psi_\varepsilon(\psi)]
			|M_\varepsilon.
}

\begin{lemma}\label{lemma:fddconv}

	For any fixed $\psi_1,\dots,\psi_n\in H^\beta_c$, consider the following Borel probability measures on $\R^n$:
$$
	\LC(\varepsilon^{-1}(f-\hat\Psi_\varepsilon)|M_\varepsilon)_n
	:=
		\LC(\varepsilon^{-1}[\langle f,\psi_1\rangle_{L^2}-\hat\Psi_\varepsilon(\psi_1),\dots,
		\langle f,\psi_n\rangle_{L^2}-\hat\Psi_\varepsilon(\psi_n)]|M_\varepsilon),
$$
where $f\sim\Pi$, and 
$$
	\mu_n:=\LC(X(\psi_1),\dots,X(\psi_n)),
$$
where $X$ is as in (\ref{X}). Then, denoting $d_{\R^n}$ the BL-metric for weak convergence on $\R^n$, we have in $P^M_{f^\dagger}$-probability as $\varepsilon\to0$.
\eq{
\label{fddconv}
	d_{\R^n}(\LC(\varepsilon^{-1}(f-\hat\Psi_\varepsilon)|M_\varepsilon)_n,\mu_n)\to 0. 
}

 \proof

	By Lemma \ref{lem:ae}, it is enough to show (\ref{fddconv}) for $f\sim\Pi^{D_\varepsilon}$, with $D_\varepsilon$ as in (\ref{newDeps}). Let $\psi\in H^\beta_c$ be fixed. Then, by taking $K=\supp(\psi)$, Lemma \ref{Laplexpansion} implies
$$
	\E^{\Pi^{D_\varepsilon}}\left[e^{\frac{\tau}{\varepsilon}[
	\langle f,\psi\rangle_{L^2}-\hat\Psi_\varepsilon(\psi)]}\big{|}M_\varepsilon\right] 
	=
		 (1+o(1)) e^{\frac{\tau^2}{2}\|L\psi\|^2_{L^2}}\frac{\Pi(D_{\varepsilon,\tau}|
		 M_\varepsilon)}{\Pi(D_{\varepsilon}|M_\varepsilon)},
$$
and the proof is concluded by showing that the ratio on the right hand side converges to 1 in $P^M_{f^\dagger}$-probability as $\varepsilon\to0$. 
Indeed, if this is the case, the convergence of the Laplace transform will imply that for any fixed $\psi\in H^\beta_c$,
\eq{
\label{singleprojconv}
	d_\R(\LC(\varepsilon^{-1}[\langle f,\psi\rangle_{L^2}-\hat\Psi_\varepsilon(\psi)]|
	M_\varepsilon),\LC(X(\psi)))\to 0,
}
in $P_{f^\dagger}^M$-probability as $\varepsilon\to0$; and, by the Cramer-Wold device, also (\ref{fddconv}) will follow by replacing $\psi$ with any finite linear combination $\sum_{i=1}^n a_i\psi_i\in H^\beta_c$.

	To proceed, first recall that Lemma \ref{lem:ae} implies that
$
\Pi(D_{\varepsilon}|M_\varepsilon)\to1
$
in $P_{f^\dagger}^M$-probability as $\varepsilon\to0$. To apply the same result to the numerator, we show that the prior probability of $D^c_{\varepsilon,\tau}$  decays exponentially as required by \eqref{eq:complementAs}. Notice
\eqstar{
	D_{\varepsilon,\tau} 
	=
		 \{f-\tau \varepsilon L\tilde\psi_\varepsilon,\ f\in G_\varepsilon\cap F_\varepsilon\}
	=
		G_{\varepsilon,\tau}\cap F_{\varepsilon,\tau},
}
where $G_{\varepsilon,\tau},\ F_{\varepsilon,\tau}$ are defined analogously to the set $D_{\varepsilon,\tau}$ introduced in the previous lemma. It is hence enough to deduce \eqref{eq:complementAs} for $G_{\varepsilon,\tau}$ and  $F_{\varepsilon,\tau}$ separately.

	First, from the definition of $G_\varepsilon$ in (\ref{setGeps}), we see that
\eqstar{
	G_{\varepsilon,\tau} 
	&\supseteq
	 	\left\{ g :\sup_{\phi\in B^\beta_K(z)}|(g,L\tilde\phi_\varepsilon \rangle_{V_\Pi}|
		\le D\sigma_\varepsilon\delta_\varepsilon/\varepsilon-
		|\tau\varepsilon| \|L\tilde\psi_\varepsilon\|_{V_\Pi}\sup_{\phi\in B^\beta_K(z)}\|			L\tilde\phi_\varepsilon \|_{V_\Pi}\right\}.
}
Now, using (\ref{sigmaepsilon}) and recalling $\delta_\varepsilon\simeq\varepsilon^{\frac{2+r-d/2}{2+r}}$,
\eqstar{
	\varepsilon\|L\tilde\psi_\varepsilon\|_{V_\Pi}
	\sup_{\phi\in B^\beta_K(z)}\|L\tilde\phi_\varepsilon \|_{V_\Pi}
	& \lesssim
		 \varepsilon^{1-\frac{\max\{0,r-\beta+4\}}{2+r}}
		 \sigma_\varepsilon
	= 
		o\left(\sigma_\varepsilon\delta_\varepsilon/\varepsilon\right).
}
Then, for all  $\varepsilon>0$ small enough 
$$
	G_{\varepsilon,\tau} 
	\supseteq 
		\left\{ g :\sup_{\phi\in B^\beta_K(z)}|(g,L\tilde\phi_\varepsilon \rangle_{V_\Pi}|
		\le D\sigma_\varepsilon\delta_\varepsilon/\varepsilon\right\},
$$
and by our particular choices of $D>\tilde D>\sqrt{6}$ we obtain (via the Borel-Sudakov-Tirelson inequality) that
$
\Pi(G^c_{\varepsilon,\tau})\le e^{-\frac{\tilde D^2}{2}(\delta_\varepsilon/\varepsilon)^2}.
$
On the other hand, for $F_\varepsilon$ defined in (\ref{Feps}),
\eqstar{
	F_{\varepsilon,\tau}
	&\supseteq 
		\{  f_1+f'_2 :  \|f_1\|_{L^2}\le \rho\varepsilon^{\frac{r-d/2}{2+r}},
		 \ \|f'_2\|_{V_\Pi}\le Q_\varepsilon+\kappa
		 -|\tau|\varepsilon\|L\tilde\psi_\varepsilon\|_{V_\Pi}\}
}
and since, by (\ref{cancx}), $\varepsilon\|L\tilde\psi_\varepsilon\|_{V_\Pi}\to0$ as $\varepsilon\to0$, we have 
$Q_\varepsilon+\kappa-|\tau|\varepsilon\|L\tilde\psi_\varepsilon\|_{V_\Pi}>Q_\varepsilon$ for $\varepsilon>0$ small enough. Then, for all such $\varepsilon>0$,
$$
	F_{\varepsilon,\tau} 
	\supseteq
		 \{  f_1+f'_2 :  \|f_1\|_{L^2}\le \rho\varepsilon^{\frac{r-d/2}{2+r}},
		  \ \|f'_2\|_{V_\Pi}\le Q_\varepsilon\},
$$
and by the isoperimetric inequality for Gaussian processes we can conclude that
$
\Pi(F^c_{\varepsilon,\tau})\lesssim e ^{-\frac{\tilde D^2}{2}\left(\delta_\varepsilon/\varepsilon\right)^2}, \ \tilde D>\sqrt{6},
$
as required.
\endproof

\end{lemma}

	Using the same argument as in the conclusion of the proof of Theorem \ref{Thm:BvM1}, we deduce from the above lemma that the semiparametric BvM phenomenon displayed in \eqref{finalconv11} does occur for all $\beta>2+d/2$, concluding the proof of Theorem \ref{thm:semiBvM2}.

\paragraph{\bf Step IV: Weak convergence in $(H^\beta_K)^*$.}

	Assume now that $\beta>2+d$. Combining the convergence of the finite dimensional distribution established in the previous step with a uniform bound on the covariance of the process $\hat X_\varepsilon$ in \eqref{hatXeps} implied by Lemma \ref{Laplexpansion}, we show that $\LC(\varepsilon^{-1}(f-\hat\Psi_\varepsilon)|M_\varepsilon)$ converges weakly to $\mu$. That $\hat\Psi_\varepsilon$ can be replaced by the posterior mean $\bar f$ can then be shown analogously as in the proof of Theorem 2.7 in \cite{Monard2019}; see Step V in the Supplement.

	It is again enough to consider the restricted prior $\Pi^{D_\varepsilon}$. Thus, for any fixed compact set $K\subset\OC$, let  $\tilde\Pi^{D_\varepsilon}(\cdot|M_\varepsilon)$ be the tight Gaussian law on $(H^\beta_K)^*$ induced by $\hat X_{\varepsilon, K}:=(\hat X_{\varepsilon}(\psi) : \psi\in H^\beta_K)$, where $\hat X_\varepsilon(\psi)$ is as in \eqref{hatXeps} with $f\sim\Pi^{D_\varepsilon}$. To exploit the convergence of the finite dimensional distributions, we further consider 'projections' of $\tilde\Pi^{D_\varepsilon}(\cdot|M_\varepsilon)$ onto suitable subspaces. In particular, let
$
	\{\Phi^\OC_{\ell r}, \ \ell\ge-1, \ r=1,\dots N_\ell\},\
	 N_\ell\lesssim 2^{\ell d},
$
be an orthonormal basis of $L^2(\OC)$ of sufficiently regular boundary corrected Daubechies wavelets. We will exploit the fact that such basis conveniently characterises the Sobolev regularity of the test functions in terms of the decay of the wavelets coefficients (see \cite{Triebel2008} or also Chapter 4 of \cite{Gine2016} for details). 

	For any $\lambda\in\N$ and all $\psi\in H^\beta_K$, let $P_\lambda\psi$ denote the projection of $\psi$ onto the finite dimensional subspace spanned by $\{\Phi^\OC_{\ell r}, \ \ell\le\lambda,\ r\le N_\ell\}$. Next, define the projected posterior $\tilde\Pi_\lambda^{D_\varepsilon}(\cdot|M_\varepsilon)$ as the law of the process
$
P_\lambda \hat X_{\varepsilon,K}:=(\hat X_\varepsilon(P_\lambda\psi) : \psi\in H^\beta_K);
$
define analogously the projected limiting law $\mu_\lambda$. For $d=d_{(H^\beta_K)^*}$ the BL-metric for weak convergence of probability measures defined on $(H^\beta_K)^*$, the triangular inequality then yields
\eq{
\label{dTriangIneq}
	d(\tilde\Pi^{D_\varepsilon}(\cdot|M_\varepsilon),\mu)
	\le
		d(\tilde\Pi^{D_\varepsilon}(\cdot|M_\varepsilon),
		\tilde\Pi_\lambda^{D_\varepsilon}(\cdot|M_\varepsilon)) 
		+ d(\tilde\Pi_\lambda^{D_\varepsilon}(\cdot|M_\varepsilon),\mu_\lambda)
		+ d(\mu_\lambda,\mu).
}

We show that the three term in the right hand side vanish. For the first, recalling the definition of the BL-metric \eqref{Eq:BLMetric}, we have
\eqstar{
	d(\tilde\Pi^{D_\varepsilon}(\cdot|M_\varepsilon),
	\tilde\Pi_\lambda^{D_\varepsilon}(\cdot|M_\varepsilon))
	&=
		\sup_{F:(H_K^\beta)^*\to\R, \ \|F\|_{Lip}\le1}\left|\E^{\Pi^{D_\varepsilon}}
		 [F(\hat X_{\varepsilon,K})-F(P_\lambda \hat X_{\varepsilon,K})]\right|\\
	&\le
		\E^{\Pi^{D_\varepsilon}}\left\|\hat X_{\varepsilon,K}-
		P_\lambda \hat X_{\varepsilon,K}\right\|_{(H^\beta_K)^*}
}
which, by definition of the norm in $(H^\beta_K)^*$, equals
\eqstar{
	\E^{\Pi^{D_\varepsilon}}
	\sup_{\psi\in B^\beta_K(1)}|\hat X_{\varepsilon}(\psi-P_\lambda \psi)|
	\le
		\E^{\Pi^{D_\varepsilon}}
		\sup_{\psi\in B^\beta_K(1)}\sum_{\ell>\lambda}\sum_{r=1}^{N_\ell} 
		|\langle \psi,\Phi^\OC_{\ell r}\rangle_{L^2}| |\hat X_{\varepsilon}(\Phi^\OC_{\ell r})|,
}
with $B^\beta_K(1)$ defined as in \eqref{Eq:Balls}. Notice that, as $\supp(\psi)\subset K$, for $\lambda$ large enough (only depending on $K$) the above sum involves only wavelets $\Phi_{\ell r}=\Phi_{\ell r}^\OC$ that are compactly supported within $\OC$. We now apply H\"older's inequality and the wavelet characterisation of Sobolev norms to upper bound the right hand side by
\eqstar{
	\E^{\Pi^{D_\varepsilon}}\sup_{\psi\in B^\beta_K(1)}
	\sum_{\ell>\lambda}\sqrt{\sum_{r=1}^{N_\ell} \langle \psi,\Phi_{\ell r}\rangle_{L^2}^2 }
	\sqrt{\sum_{r=1}^{N_\ell}|\hat X_{\varepsilon}(\Phi_{\ell r})|^2}
	\lesssim
		\sum_{\ell>\lambda}2^{-\ell\beta} \E^{\Pi^{D_\varepsilon}}
		\sqrt{\sum_{r=1}^{N_\ell}|\hat X_{\varepsilon}(\Phi_{\ell r})|^2}.
}
Jensen's inequality implies the further upper bound
\eqstar{
	\sum_{\ell>\lambda}2^{-(\beta-\beta') \ell}
	\sqrt{\sum_{r=1}^{N_\ell}\E^{\Pi^{D_\varepsilon}}
	|\hat X_{\varepsilon}(2^{-\beta'\ell}\Phi_{\ell r})|^2}
}
having scaled the wavelets by a factor $2^{-\beta'\ell}$ for some $2+d/<\beta'<\beta-d/2$ (possible since we here assume $\beta>2+d$).
In particular, since
$
	\|\Phi_{\ell r}\|^2_{H^{\beta'}}\simeq 2^{2\beta' \ell},
$
we have that $2^{-\beta' \ell}\Phi_{\ell r}\in B^{\beta'}_c(1)$. Using Lemma \ref{Laplexpansion} and the fact that $\|L\Phi_{\ell r}\|^2_{L^2}\lesssim\|\Phi_{\ell r}\|_{H^2}\simeq 2^{4l}$, we obtain
\eqstar{
	\E^{\Pi^{D_\varepsilon}}
	\left[e^{\hat X_\varepsilon(2^{-\beta' \ell}\Phi_{\ell r})}\big{|}M_\varepsilon\right]
	&=
		e^{R_\varepsilon}e^{\frac{1}{2}\|2^{-\beta' \ell} L \Phi_{\ell r}\|^2_{L^2}}
		\frac{\Pi(D_{\varepsilon,\tau}|M_\varepsilon)}{\Pi(D_\varepsilon|M_\varepsilon)}\\
	&\le
		\frac{e^{R_\varepsilon}}{\Pi(D_\varepsilon|M_\varepsilon)}
		e^{2^{-2\beta' \ell-1} \|L \Phi_{\ell r}\|^2_{L^2}}\\
	&\lesssim
		 r_\varepsilon\\
	&=
		O_{P^M_{f^\dagger}}(1)
}
having lower bounded the probability on the denominator by Lemma \ref{lem:ae}. Then, since $x^2\le e^x+e^{-x}$ for all $ x\in\R$,
$$
	\E^{\Pi^{D_\varepsilon}} |\hat X_{\varepsilon}(2^{-\beta' \ell}\Phi_{\ell r})|^2
	\le
		 \E^{\Pi^{D_\varepsilon}} e^{\hat X_{\varepsilon}(2^{-\beta' \ell}\Phi_{\ell r})}
		 +\E^{\Pi^{D_\varepsilon}} e^{-\hat X_{\varepsilon}(2^{-\beta' \ell}\Phi_{\ell r})}
	\lesssim 
		r_\varepsilon.
$$
We thus obtain, recalling $N_\ell\simeq 2^{\ell d},$
\eqstar{
	d(\tilde\Pi^{D_\varepsilon}(\cdot|M_\varepsilon),
	\tilde\Pi_\lambda^{D_\varepsilon}(\cdot|M_\varepsilon))
	\lesssim
		\sum_{\ell>\lambda}2^{-(\beta-\beta') \ell}\sqrt{\sum_{r=1}^{N_\ell}r_\varepsilon}
	\lesssim 
		r'_\varepsilon \sum_{\ell>\lambda} 2^{-(\beta-\beta'-d/2)\ell}.
}
Since we have chosen $\beta'<\beta-d/2$, the latter series is convergent, implying that the right hand side vanishes if $\lambda\to\infty$.

For the second term in \eqref{dTriangIneq}, we can deduce directly from Lemma \ref{lemma:fddconv} that for any fixed $\lambda$ we have 
$
d(\tilde\Pi_\lambda^{D_\varepsilon}(\cdot|M_\varepsilon),\mu_\lambda)\to0
$
in $P^M_{f^\dagger}$-probability as $\varepsilon\to0$. For the third term we proceed similarly to the first, obtaining
\eqstar{
	d(\mu,\mu_\lambda)
	=
		E\sup_{\psi\in B^\beta_K(1)}| X(\psi-P_\lambda \psi)|
	\lesssim
		\sum_{\ell>\lambda}2^{-\beta \ell}\sqrt{\sum_{r=1}^{N_\ell}E |X(\Phi_{\ell r})|^2}.
}
Recall $X(\Phi_{\ell r})\sim \NC(0,\|L\Phi_{\ell r}\|^2_{L^2})$. Hence
\eqstar{
	d(\mu,\mu_\lambda)
	\lesssim
	\sum_{\ell>\lambda}2^{-\beta \ell}2^{2l}\sqrt{N_\ell}
	\le
	\sum_{\ell>\lambda}2^{-(\beta-2-d/2) \ell}
}
which again converges since $\beta>2+d$.  To conclude, we can fix arbitrarily $\varepsilon'>0$, and then find $\lambda=\lambda(\varepsilon')$ sufficiently large so that the first and third term  in (\ref{dTriangIneq}) are smaller than $\varepsilon'$. For such value of $\lambda$, the second term can be made smaller that $\varepsilon'$ by choosing $\varepsilon$ small enough with $P^M_{f^\dagger}$-probability approaching one, implying that
$
	d(\tilde\Pi^{D_\varepsilon}(\cdot|M_\varepsilon),\mu)\to0 
$
in $P_{f^\dagger}^M$-probability as $\varepsilon\to0$.

\appendix

\section{Information lower bound for linear inverse problems}
\label{appendix:infolb}

	Let $M_\varepsilon$ be given by (\ref{eq:observation}) with $f=f^\dagger$, and let $\ell(f)=\log p_{f}(M_\varepsilon), \ f\in \overline W,$ be the log-likelihood in (\ref{Likelihood}). For any $ h\in \Wt,\ \varepsilon>0$, we have
\eqstar{
	\log\frac{p_{f^\dagger+\varepsilon h}(M_\varepsilon)}{p_{f^\dagger}(M_\varepsilon)} 
	= 
		\ell(f^\dagger+\varepsilon h)-\ell(f^\dagger)
	=
		\langle A h,\W\rangle_{W_2}-\frac{1}{2}\| Ah\|^2_{W_2}.
}
Recalling $\langle A h,\W\rangle_{W_2}\sim\NC(0,\| Ah\|^2_{W_2})$, the model is seen to be locally (asymptotically) normal (LAN), with LAN-inner product and norm respectively given by
$$
	\langle\cdot,\cdot\rangle_{LAN}
	=
		\langle A\cdot,A\cdot\rangle_{W_2}, \quad \|\cdot\|_{LAN}=\|A\cdot\|_{W_2}.
$$

	Let $\psi\in W_1$ satisfy the assumptions of Theorem \ref{Thm:BvM1}, and consider the continuous linear map
$$
	L_\psi: \Wt\to \R, \quad  L_\psi(h)
	=
		\langle h,\psi\rangle_{W_1},
$$
defined by extension using the fact that $W_1\subseteq\overline W$ is dense.
As by assumption $\psi=-A^*A\tilde\psi$ for some $\tilde\psi\in V_\Pi$, then for all $h\in \Wt$
\eqstar{
	L_\psi(h)
	&=
		\langle \psi,h\rangle_{W_1}
	=
		\langle -A^*A\tilde\psi,h\rangle_{W^2}
	= 
		\langle-\tilde\psi,h\rangle_{LAN},
}
so that the Riesz representer with respect to the LAN-inner product of the linear functional $L_\psi$ is $-\tilde\psi$. We then deduce from the semiparametric theory of efficiency (see Chapter 25 in \cite {Vaart1998}, or Section 7.5 in \cite{Nickl2017}) that the information lower bound for estimating $L_\psi(f^\dagger)=\langle f^\dagger,\psi\rangle_{W_1}$ from model (\ref{eq:observation}) is identified by the random variable
$$
	Z\sim\NC(0,\|\tilde\psi\|_{LAN})=\NC0,\| A\tilde\psi\|_{W_1}), 
$$
and we have the lower bound \eqref{infolb1}  for the asymptotic minimal variance . Note that when $A^*A$ has a well defined inverse we can write $\|A\tilde\psi\|^2_{W_2}=\|(A^*A)^{-1}\psi\|^2_{LAN}$. In analogy with the finite-dimensional case, we then sometimes call $A^*A$ the Fisher information operator.

\section{Properties of elliptic boundary value problems}
\label{Sec:EllipticBVPFacts}

We list here some key facts relative to the BVP (\ref{BVP}) following from the general elliptic theory (see, e.g., the monographs \cite{Lions1972,Roitberg2012}). We start noting that the operator $L$ defines a bounded linear operator from $H^s=H^s(\OC)$ into $H^{s-2}$ for all $s\ge 2$,
 and, in view of the symmetry of the coefficients $a_{ij}$, it is also self-adjoint with respect to $\langle \cdot,\cdot\rangle_{L^2}$ when acting upon $H^2_0$. 
 
 	If, in addition, we assume the uniform ellipticity condition:
$$
	\sum_{i,j=1}^da_{ij}(x)\xi_i\xi_j\ge\theta|\xi|^2,
	\quad
	 \forall x\in\OC,\ \xi=(\xi_1,\dots,\xi_d)\in\R^d,
$$
for some constant $\theta>0$, then for all $s\ge0$ and  any $f\in H^s$, there exists a unique weak solution $u_f\in H_0^{s+2}$ of (\ref{BVP}) satisfying the variational formulation of the problem:
\eq{\label{varfor}
	\int_\OC\sum_{i,j=1}^da_{ij}\frac{\partial u_f}{\partial x_i} \frac{\partial v}{\partial x_j}
	=
		\int_\OC fv, \quad \forall v\in H^1_0.
}
Furthermore, we have the elliptic estimates
$$
	\|u_f\|_{H^{s+2}}\le c_s \|f\|_{H^s},
$$
for constants $s_s>0$ depending only on $s$. These results follow directly from Theorem 5.4 in \cite[Chapter 2]{Lions1972} (see also remark 7.2 in the same reference) by noting that $u_f=0$ is the unique smooth solution of \eqref{BVP} with $f=0$ (e.g., in view of Theorem 3 and Theorem 4 in \cite[Section 6.3]{Evans1998}).
Finally, as pointed out in Remark (ii), page 310 in \cite{Evans1998},  it follows that $Lu_f=f$ almost everywhere on $\OC$.

	With a slight abuse of notation, let $L^{-1}$ denote the solution map, so that $L^{-1}f=u_f$ is the unique element in $H^{s+2}_0$ that satisfies (\ref{varfor}). Also, in view of the uniqueness of weak solutions, $L^{-1}Lu=u$ for all $u\in H^2_0$. From the above results we have that 
$
	L^{-1}:H^s\to H_0^{s+2},\ s\ge0  
$
defines a linear and bounded isomorphism which is self-adjoint on $L^2$ (following from the self-adjointness of $L$).

\paragraph{Acknowledgements.}
The authors would like to thank Richard Nickl for suggesting the problem and for his valuable assistance. M.G. also wishes to thank Sven Wang for helpful discussions. The authors are grateful to the Associate Editor and two referees for helpful comments that have greatly improved the paper. This work was supported by the European Research Council under ERC grant UQMSI (No. 647812). In addition M.G. was supported by the UK Engineering and Physical Sciences Research Council (EPSRC) grant EP/L016516/1 for the Cambridge Centre for Analysis, and H.K. was supported  by the Cantab Capital Institute for the Mathematics of Information (RG83535).

%
%
%
%
%

{\small
\bibliography{Inverse_problems_references}}
\bibliographystyle{plain}

%
%
%
%
%

\section{Supplementary materials}
\label{App:RemainingProofs}

	In this supplement, we provide the remaining proofs for the main paper and prove the auxiliary results used to show the general semiparametric Bernstein-von Mises theorem.

\section*{1 $\ \ \ $ Remaining proofs}
\label{App:RemainingProofs}

%
%
%

\subsection*{1.1 $\ \ $ Proof of Corollary \ref{Cor:SemiUQ}.}
\label{ProofCorUQ}

	The proof follows the argument in Section 2 of \cite{Castillo2013}. We start by noting that the function $\Phi:[0,\infty)\to[0,1]$ is uniformly continuous and strictly increasing, with continuous inverse $\Phi^{-1}:[0,1]\to[0,\infty)$. Thus, for every $\gamma>0$ we can find $\varepsilon>0$ such that
$
	|\Phi(t+\varepsilon)-\Phi(t)|\le \gamma,\ \forall t\ge0.
$
For such $\varepsilon$ and for all $t\ge0$, denoting $B(0,t)=\{x\in\R, \ |x|\le t\}$, we have
\eqstar{
	\Pr (t-\varepsilon<|Z|\le t+\varepsilon)
	=
		 \Phi(t+\varepsilon)-\Phi(t-\varepsilon)
	\le
		 2\gamma.
}
Thus, applying Corollary 7.3.22 in \cite{Gine2016} to 
$\LC(\varepsilon^{-1}\langle f-\bar f,\psi\rangle_{W_1}|M_\varepsilon)$ converging weakly to $\LC(Z)$ in $P^M_{f^\dagger}$-probability as $\varepsilon\to0$, we deduce that
\eqstar{
	\sup_{0\le t<\infty}
	&|\Pi(\varepsilon^{-1}\langle f-\bar f,\psi\rangle_{W_1}\in B(0,t)\g M_\varepsilon)
	-\Pr(Z\in B(0,t))|
	=
		o_{P^M_{f^\dagger}}(1)
}
as $\varepsilon\to0$. Thus, recalling the definition of $R_\varepsilon$ after (\ref{Ceps})
\eqstar{
	\Phi(\varepsilon^{-1}R_\varepsilon)   
	&=
		 \Pr(|Z|\le \varepsilon^{-1}R_\varepsilon)- \Pi(\varepsilon^{-1}|\langle f-\bar f,
		 \psi\rangle_{W_1}|\le \varepsilon^{-1} R_\varepsilon \g M_\varepsilon) 
		 + 1-\alpha\\
	&=
		1-\alpha+o_{P^M_{f^\dagger}}(1)
}
as $\varepsilon\to0$ by the above with $t=\varepsilon^{-1}R_\varepsilon$. Since $\Phi^{-1}$ is continuous, the continuous mapping theorem yields
$$
	\varepsilon^{-1}R_\varepsilon
	=
		\Phi^{-1}[\Phi(\varepsilon^{-1}R_\varepsilon)]\to^{P^M_{f^\dagger}}\Phi^{-1}(1-\alpha).
$$
Then the first claim follows using Theorem \ref{Thm:BvM2}, as
\eqstar{
	P^M_{f^\dagger}(\langle f^\dagger,\psi\rangle_{W_1}\in C_\varepsilon)
	&= 
		P^M_{f^\dagger}(\varepsilon^{-1}|\langle f^\dagger-\bar f,\psi\rangle_{W_1}|
		\le \varepsilon^{-1}R_\varepsilon )\\
	&=
		 P^M_{f^\dagger}\left(\varepsilon^{-1}|\langle f^\dagger-\bar f,\psi\rangle_{W_1}|
		 \le \Phi^{-1}(1-\alpha)\right)+o(1)\\
	&=
		 \Pr\left(|Z|\le \Phi^{-1}(1-\alpha)\right)+o(1)\\
	&=
		1-\alpha+o(1).
}

%
%
%

\subsection*{1.2 $\ \ $ Proof of Lemma \ref{lemma:tightness}.}
\label{ProofTightness}

	1. First assume that $\beta>2+d/2$, and denote $B^\beta_c(1):=\{h\in H^\beta_c,$ $\|h\|_{H^\beta}\le1\}$. According to (\ref{X}), $X_{|B^\beta_c(1)}:=(X(\psi) : \psi\in B^\beta_c(1))$ is a Gaussian process with intrinsic distance
\eqstar{
	d^2_X(\psi,\psi')
	&:=
		E[X(\psi)-X(\psi')]^2
	\lesssim
		\|\psi-\psi'\|^2_{H^2}.
}
Next, from Edmund and Triebel's upper bound for the entropy numbers in general Besov spaces 
(see \cite{Triebel1978})  we deduce that, for positive reals $s_1<s_2$, denoting $B^s(r):=\{h\in H^s, \ \|h\|_{H^s}\le r\}, \ r>0,$
\eq{
\label{Edbound1}  
	\log N(\eta,B^{s_2}(1),\|\cdot\|_{H^{s_1}})
	\lesssim
		 \eta^{-\frac{d}{s_2-s_1}},
	\quad 
		\eta>0.
}
Then it follows from Dudley's metric entropy inequality \cite[Theorem 2.3.7]{Gine2016} that for all $z>0$
\eq{\label{metricbound}
	E\sup_{\psi\in B_c^\beta(1),\ \|\psi\|_{H^2}\le z}|X(\psi)|
	 \lesssim
	 	 \int_0^z\sqrt{2\log N(\eta, B^\beta(1),\|\cdot\|_{H^2}})d\eta
 	\lesssim
		\int_0^z \eta^{-\frac{d}{2(\beta-2)}} d\eta,
}
which is indeed convergent for all $\beta>2+d/2$. Thus, letting $z\to0$ in (\ref{metricbound}) implies that $X_{|B^\beta_c(1)}$ has a version taking values in the separable 
Banach space   
\eq{\label{subspaceB}
(\B,\|\cdot\|_\B),\ \B=UC(B^\beta_c(1),d_X), \quad \|x\|_\B=\sup_{\psi\in B^\beta_c(1)}|x(\psi)|,
}
of bounded and uniformly continuous (with respect to the metric $d_X$ on $B^\beta_c(1)$) pre-linear functionals on $B^\beta_c(1)$, the separability following from Corollary 11.2.5 in \cite{Dudley2002} since, in view of (\ref{Edbound1}),  $B^\beta_c(1)$ is totally bounded for the metric $d_X$ if $\beta>2$. Finally, as according to (\ref{subspaceB}) $\B$ is an isometrically imbedded closed subspace of $(H^\beta_c)^*$, we deduce from Oxtoby-Ulam theorem (Proposition 2.1.4 in \cite{Gine2016}) that $X_{|B^\beta_c(1)}$ induces a tight Borel Gaussian probability measure on $\B$, which has a unique  extension to $(H^\beta_c)^*$.

	2. For $\beta<2+d/2$, as  $H^\beta_c\subset H^{\beta'}_c$ with continuous embedding if $\beta'<\beta$, it is enough to show that
$
	\Pr\left(\sup_{\psi\in B^\beta_c(1)}|X(\psi)|<\infty \right)=0
$
for $2<\beta<2+d/2$. We proceed by contradiction, assuming on the contrary that
\eq{\label{contr1}
	\Pr\left(\sup_{\psi\in B_c^\beta(1)}|X(\psi)|<\infty\right)>0.
}
In view of (\ref{Edbound1}), $B_c^\beta(1)$ is separable with respect to the intrinsic metric $d_X$ for any $\beta>2$. Hence Proposition 2.1.12 in \cite{Gine2016} and (\ref{contr1}) jointly imply, by Proposition 2.1.20 in \cite{Gine2016}, that
$
	E\sup_{\psi\in B_c^\beta(1)}|X(\psi)|<\infty,
$
which we will show to yield a contradiction. To do so, notice that $(X(L^{-1}\psi)$ : $ \psi\in H^\beta_c)$ has the same law on $\R^{H^\beta_c}$ as the standard gaussian white noise $\W$. 
Thus,
$$
	E\sup_{\psi\in B_c^\beta(1)}|X(\psi)|
	=
		 \E\sup_{\psi\in B_c^\beta(1)}|\langle L\psi,\W\rangle_{L^2}|,
$$
and the proof is completed by  finding a suitable lower bound to show that the right hand side diverges.

	Considering the orthonormal Daubechies wavelet basis of $L^2$ introduced in Step IV in Section \ref{sec:MtheoProof},  select for each $j\ge1$, $n_j=c'2^{jd}, \ c'>0,$ wavelets $\{\Phi^\OC_{jr}, \ r=1,\dots,n_j\}$ with disjoint compact support within $\OC$. Next, for each $m=1,\dots,2^{n_j}$ and $b_{m\cdot}=(b_{mr},  \ r=1,\dots,n_j) \in\{-1,1\}^{n_j}$, define
\eq{\label{hmj}
	h_{jm}(x)=k_j\sum_{r=1}^{n_j}b_{mr}2^{-j\beta}\Phi^\OC_{jr}(x), \quad x\in\OC,
}
where $k_j>0$ is to be fixed. Recall that it is enough to consider $2<\beta<2+d/2$. We have $h_{jm}\in  H_c^\beta$, and by the usual wavelet characterisation  of the Sobolev norms
\eqstar{  
	\|h_{jm}\|^2_{H^\beta}
	\simeq
		 \sum_{\ell\ge-1}\sum_{s=1}^{n_j}2^{2l\beta}
		 \langle h_{jm},\Phi^\OC_{ls}\rangle^2_{L^2}
	=
		k^2_jn_j.
}
Hence, choosing $k_j<n_j^{-1/2}$ guarantees that 
$\{h_{jm},\ m=1,\dots,2^{n_j}\}$ $\subset B^\beta_c(1)$, yielding the lower bound
$$
	\E\sup_{\psi\in B_c^\beta(1)}|\langle L\psi,\W\rangle_{L^2}| 
	\ge
		 \E\max_{m=1,\dots,2^{n_j}}|\langle Lh_{mj},\W\rangle_{L^2}|,
	\quad
		 j\ge1,
$$
which we can further develop by restricting the maximum to a suitable smaller subset. In particular, the Gaussian vector $(\W(Lh_{jm}), \ m=1,\dots,2^{n_j})$ has intrinsic metric 
\eqstar{
	d_j^2(h_{jm},h_{jm'}) 
	&=
		 \|L(h_{jm}-h_{jm'})\|^2_{L^2}
	=
		 k^2_j 2^{-2j\beta}\left\| \sum_{r=1}^{n_j} (b_{mr}-b_{m'r}) 
		 L\Phi^\OC_{jr} \right\|^2_{L^2};
}
and arguing as in the proof of Proposition 6 in \cite{Nickl2017} we can select, for sufficiently large $j$,  a subset $\{h_{j1},\dots,h_{jm_j}\}\subseteq\{h_{jm},$ $ m=$ $1,\dots,$ $2^{n_j}\}$ of cardinality $m_j\ge 3^{n_j/4}$, such that
$$
	d^2_j(h_{jh},h_{jk}) \gtrsim  2^{2j(2-\beta)}, \quad h\neq k.
$$
Thus, by applying Sudakov's lower bound \cite[Theorem 2.4.12]{Gine2016}, we deduce that for all such $j$
\eqstar{
	\E\max_{m=1,\dots,2^{n_j}}|\langle Lh_{jm},\W\rangle_{L^2}| 
	&\ge
		 \E\max_{h=1,\dots,m_j}\{|\langle Lh_{jh},\W\rangle_{L^2}|\} \\
	&\ge
		 c 2^{j(2-\beta)}\sqrt{\log N(2^{j(2-\beta)},\{h_{j1},\dots,h_{jm_j}\},d_j)}\\
	&\ge
		 c' 2^{j(2-\beta)}\sqrt{\log m_j}\\
	&\ge
		 c''2^{j(2+d/2-\beta)}. 
}
The last line diverges as $j\to\infty$ for all $\beta<2+d/2$, yielding the contradiction.

	3. Assuming tightness on $(H^\beta_c)^*$ for $\beta=2+d/2$ would imply (exactly as above) that $X$ were sample bounded and, in view of Proposition 2.1.7 in \cite{Gine2016}, also sample continuous with respect to $d_X$. Then, Proposition 2.4.14 in \cite{Gine2016} would yield 
$$
\lim_{\eta\to0}\eta\sqrt{\log N(\eta,B^\beta_c(1),d)}=0
$$
which, taking the sequence $\eta_j=2^{j(2-\beta)}=2^{-jd/2}$, is in contrast with the fact that
\eqstar{
2^{j(2-\beta)}\sqrt{\log N(2^{j(2-\beta)},B^\beta_c(1),d)} 
\ge 2^{j(2-\beta)}\sqrt{\log N(2^{j(2-\beta)},\{h_{j1},\dots,h_{jm_j}\},d_j)},
}
and that the right hand side is bounded below by a positive constant for $\beta=2+d/2$, as seen above.

\qed    

%
%
%
%
%

\section*{2 $\ \ \ $ Proof of supporting lemmas for Theorem \ref{Thm:BvM1}}
\label{App:ProofAuxiliaryResults}

%
%
%

\subsection*{2.1 $\ \ $ Proof of Lemma \ref{lem:ae}.}
\label{App:ProofLemAe}

	We start by noting that $\Pi(B)=\Pi(B\cap D_\varepsilon)+\Pi(B\cap D_\varepsilon^c)$ and
\begin{align*}
	\Pi(B\cap D_\varepsilon)-\Pi^{D_\varepsilon}(B) 
	=
		 \frac{\Pi(B\cap D_\varepsilon)}{\Pi(\Wt)}-\frac{\Pi(B\cap D_\varepsilon)}				{\Pi(D_\varepsilon)}=-\Pi(D_\varepsilon^c)\Pi^{D_\varepsilon}(B)
\end{align*}
which implies $\|\Pi(\cdot \g M_\varepsilon)-\Pi^{D_\varepsilon}(\cdot \g M_\varepsilon)\|_{TV}\leq 2\Pi(D_\varepsilon^c\g M_\varepsilon)$.  Hence it suffices to prove the first limit in \eqref{Eq:Theo1Thesis2}. This will be done using Markov's inequality and showing that $E^M_{f^\dagger}(\Pi(C^c_\varepsilon\g M_\varepsilon))\to0$. In particular, we split the expectation as
\begin{align}
\label{eq:ConvergenceExpectation}
	E_{f^\dagger}^M(\Pi(D_\varepsilon^c  |  M_\varepsilon)) 
	=
		 E_{f^\dagger}^M(\Pi(D_\varepsilon^c  |  M_\varepsilon)
		 1_{F_\varepsilon})
		 +E_{f^\dagger}^M(\Pi(D_\varepsilon^c  |  M_\varepsilon)
		  1_{F_\varepsilon^c})
\end{align}
where $F_\varepsilon$ is a suitable event to be specified for which $P^M_{f^\dagger}(F_\varepsilon)\to0$, yielding the cancellation of the first term, at a sufficiently slow rate so that also the second vanishes due to the assumption on $\Pi(D^c_\varepsilon)$.

	We proceed constructing $F_\varepsilon$. For $\ell(f)$ the log-likelihood defined in \eqref{Likelihood}, we can rewrite the posterior  \eqref{eq:BayesFormula} as
\begin{align}
\label{Eq:Post2}
	\Pi(B |  M_\varepsilon) 
	& =
		 \frac{\int_B e^{\ell(f)-\ell(f^\dagger)}d\Pi(f)}{\int_{\Wt} e^{\ell(f)-\ell(f^\dagger)}d\Pi(f)}
		  \quad B\in\mathcal{B}_{\Wt}. 
\end{align}  
It follows from \eqref{Likelihood} that under $P^M_{f^\dagger}$ we have 
\begin{align*}
\ell(f)-\ell(f^\dagger) & = \frac{1}{\varepsilon}\langle A(f-f^\dagger),\W\rangle _{W_2}- \frac{1}{2\varepsilon^2}\|A(f-f^\dagger)\|^2_{W_2}. 
\end{align*}

	Let $\nu$ be any probability measure on the set $B=\{ f : \|A(f-f^\dagger)\|\leq\delta_\varepsilon\}$. Applying Jensen's inequality to the exponential function we get for any $\tilde C>-1/2$
\begin{align*}
	P_{f^\dagger}^M&\left(\int_B e^{\ell(f)-\ell(f^\dagger)}d\nu(f)\leq e^{-(1+\widetilde{C})
	\left(\delta_\varepsilon/\varepsilon\right)^2}\right)\\
	& \leq
		 \mathbb{P} \left(\E^\nu\left( \frac{1}{\varepsilon}\langle A(f-f^\dagger),
		 \W\rangle _{W_2} - \frac{1}{2\varepsilon^2}\|A(f-f^\dagger)\|^2_{W_2}\right)
		 \leq {-(1+\widetilde{C})\left(\delta_\varepsilon/\varepsilon\right)^2}\right). 
\end{align*}
Denote $Z= \frac{1}{\varepsilon}\int_B \langle A(f-f^\dagger),\W\rangle _{W_2} d\nu(f)\sim \NC(0,C_Z)$ where, using again Jensen's inequality,
\begin{align*}
	C_Z 
	& =\frac{1}{\varepsilon^2} \E\left(\E^\nu \langle A(f-f^\dagger),
		 \W\rangle _{W_2} \right)^2\\
	& \leq 
		\frac{1}{\varepsilon^2} \E^\nu \left(\E \langle A(f-f^\dagger),
		\W\rangle _{W_2}^2\right) \\
	&  =
		\frac{1}{\varepsilon^2} \int_B \| A(f-f^\dagger)\| _{W_2}^2 d\nu(f)\\
	&\leq
		\left(\delta_\varepsilon/\varepsilon\right)^2. 
\end{align*}
We can then conclude 
\begin{align*}
	P_{f^\dagger}^M\left(\int_B e^{\ell(f)-\ell(f^\dagger)}d\nu(f)
	\leq 
		e^{-(1+\widetilde{C})\left(\delta_\varepsilon/\varepsilon\right)^2}\right)
	 &=
		\mathbb{P}\left(\left| Z-\E Z \right|  \ge 
		\Big(\frac{1}{2}+\widetilde{C}\right)\left(\delta_\varepsilon/\varepsilon\right)^2\Big)\\	&\leq
		 e^{-\frac{(1/2+\tilde{C})^2}{2}(\delta_\varepsilon/\varepsilon)^2}
\end{align*}
the last inequality following from the standard Gaussian tail bound $\P(|Z-\E Z| \ge c)\leq e^{-c^2/(2\mathbb{V}\text{ar}(Z))}$. We can now choose $\nu=\Pi(\cdot\cap B)/\Pi(B)$ and let 
\begin{align*}
	F_\varepsilon
	=
		\left\{f:\int_B e^{\ell(f)-\ell(f^\dagger)}d\nu(f)
		\leq e^{-\frac{3}{2}(\delta_\varepsilon/\varepsilon)^2}\right\}.
\end{align*}  
Using the above with $\tilde C=1/2$ we see that $P_{f^\dagger}^M(F_\varepsilon)\leq e^{-\frac{1}{2}(\delta_\varepsilon/\varepsilon)^2}\to0$, which implies that the first term in \eqref{eq:ConvergenceExpectation} tends to zero since $\Pi(\cdot  |  M_\varepsilon)\leq1$. 

%

	For the second term study the small ball probabilities $\Pi(B)=\Pi(f : \|A(f-f^\dagger)\|_{W_2}\leq\delta_\varepsilon)$ using the condition \eqref{eq:concentrationAs} on the concentration function. We see from \eqref{Eq:Post2}
\begin{align*}
	E_{f^\dagger}^M(\Pi(D_\varepsilon^c  |  M_\varepsilon)1_{F_\varepsilon^c})  
	& \leq 
		E_{f^\dagger}^M\left(\frac{\int_{D_\varepsilon^c}e^{\ell(f)-\ell(f^\dagger)}d\Pi(f)}
		{\int_{B}e^{\ell(f)-\ell(f^\dagger)}\Pi(B)d\nu(f)}\,1_{F_\varepsilon^c}\right)\\     
	& \leq
		 \frac{e^{2(\delta_\varepsilon/\varepsilon)^2}}{\Pi(f : \|A(f-f^\dagger)\|_{W_2}^2
		 \leq\delta_\varepsilon^2)}\int_{D_\varepsilon^c}E_{f^\dagger}^M
		 \left( e^{\ell(f)-\ell(f^\dagger)}\right) d\Pi(f).
\end{align*}
For $f\sim\Pi$, denote the concentration function of the (image) Gaussian measure $\tilde\Pi=\LC(Af)$ as
\begin{align*}
	\widetilde{\phi}_{\Pi,f^\dagger}(\delta)
	=
		\inf_{Ag\in V_{\widetilde{\Pi}},\ \|A(g-f^\dagger)\|_{W_2}\leq\delta}
		\frac{\|Ag\|_{V_{\widetilde{\Pi}}}^2}{2}-\log\Pi(f : \|Af\|_{W_2}\leq\delta).
\end{align*}
Following Proposition 2.6.19 in \cite{Gine2016} we next show that   
\begin{align*}
	\Pi(f : \|A(f-f^\dagger)\|_{W_2}^2\leq\delta^2) 
	\geq 
		e^{-\widetilde{\phi}_{\Pi,f^\dagger}(\delta/2)}. 
\end{align*}
Let $g\in V_\Pi$ be such that $\|A(g-f^\dagger)\|_{W_2}\leq \delta/2$. Then $\|A(f-f^\dagger)\|_{W_2}\leq \|A(f-g)\|_{W_2}+ \delta/2$. We denote $\Pi_g(B)=\Pi(B-g)=\Pi(f: f+g\in B)$.   
Using the Cameron-Martin theorem \cite[Corollary 2.4.3.]{Bogachev1998}  and the fact that $f$ is a centred Gaussian random variable we can write  
\begin{align*} 
	\Pi(&f : \|A(f-f^\dagger)\|_{W_2}\leq\delta)
 	\geq 
		\Pi(f : \|A(f-g)\|_{W_2}\leq \delta/2)\\  
	& =
		 \frac{1}{2}\left(\Pi_{-Ag}\left(f : \|Af\|_{W_2}\leq \frac{\delta}{2}\right)
		 +\Pi_{Ag}\left(f : \|Af\|_{W_2}\leq \frac{\delta}{2}\right) \right)\\
	& =
		 \frac{1}{2}\left(\int_{\{\|\widetilde{f}\|_{W_2}\leq \frac{\delta}{2}\}}
		 \frac{d\widetilde{\Pi}_{-\widetilde{g}}(\widetilde{f})}{d\widetilde{\Pi}(\widetilde{f})}
		 d\widetilde{\Pi}(\widetilde{f})
		+\int_{\{\|\widetilde{f}\|_{W_2}\leq \frac{\delta}{2}\}}\frac{d\widetilde{\Pi}
		_{\widetilde{g}}(\widetilde{f})}{d\widetilde{\Pi}(\widetilde{f})}
		d\widetilde{\Pi}(\widetilde{f})\right)\\
	& = 
		\frac{1}{2}\left(\int_{\{\|\widetilde{f}\|_{W_2}\leq \frac{\delta}{2}\}} 
		\left(e^{-\langle \widetilde{g},\widetilde{f}\rangle_{V_{\widetilde{\Pi}}}}
		+e^{\langle \widetilde{g},\widetilde{f}\rangle_{V_{\widetilde{\Pi}}}}\right)
		e^{-\frac{\|\widetilde{g}\|_{V_{\widetilde{\Pi}}}^2}{2}}d\widetilde{\Pi}
		(\widetilde{f})\right)\\
	& \geq 
		e^{-\frac{\|Ag\|_{V_{\widetilde{\Pi}}}^2}{2}} 
		\Pi\left(f : \|Af\|_{W_2}\leq\frac{\delta}{2}\right)
\end{align*}
where $Ag=\widetilde{g}$ and $Af=\widetilde{f}$. The last inequality follows from the fact $e^{-x}+e^x\geq 2$ for all $x\in\R$. 
We can then conclude   
\begin{align*}
	E_{f^\dagger}^M(\Pi(D_\varepsilon^c \g M_\varepsilon)1_{F_\varepsilon^c})
	\leq 
		e^{2(\delta_\varepsilon/\varepsilon)^2}
		e^{\widetilde{\phi}_{\Pi,f^\dagger}(\delta_\varepsilon/2)}\Pi(D_\varepsilon^c)
\end{align*}
since $E_{f^\dagger}^M\left( e^{\ell(f)-\ell(f^\dagger)}\right)=1$.

	Note that $A$ is assumed to be linear and injective and hence the RKHS $V_{\widetilde{\Pi}}=A(V_\Pi)$ of $Af$ is isometric to $V_\Pi$ (see Exercise 2.6.5 in \cite{Gine2016}). By assumption,
$
\|Af\|_{W_2}\leq c\|f\|_{\Wt}
$
for all $f\in\Wt$, which implies 
\begin{align*}
	-\log\Pi(f  : \|Af\|_{W_2}\leq \delta)
	\leq 
		-\log\Pi\left(f : \|f\|_{\Wt}\leq \delta/c\right). 
\end{align*}
We also have $\|A(g-f^\dagger)\|_{W_2}\leq c\|g-f^\dagger\|_\Wt$ and hence $\widetilde{\phi}_{\Pi,f^\dagger}(\delta)\leq\phi_{\Pi,f^\dagger}(\delta/c)$ for all $\delta$. Thus by \eqref{eq:concentrationAs} and assumption \eqref{eq:complementAs} we can conclude
\begin{align*}
	E_{f^\dagger}^M(\Pi(D_\varepsilon^c \g M_\varepsilon)1_{F_\varepsilon^c})
	\leq
		 e^{2(\delta_\varepsilon/\varepsilon)^2}
		 e^{{\phi}_{\Pi,f^\dagger}(\delta_\varepsilon/2c)}\Pi(D_\varepsilon^c)
		 \leq e^{(3-D)(\delta_\varepsilon/\varepsilon)^2}
	\to
		0.
\end{align*}

\qed

%
%
%

\subsection*{2.2 $\ \ $ Proof of Lemma \ref{lem:Convergence}.}
\label{App:ProofLemConv}

	Denote $f_{\tau}=f+\tau\varepsilon\pt$. Then the left hand side of \eqref{eq:expected} can be written as 
\begin{align*}   
	\E^{\Pi^{D_\varepsilon}}&\left(e^{\frac{\tau}{\varepsilon}\langle f-f^\dagger,
	\psi\rangle_{W_1}+\tau\langle A\pt,\W\rangle_{W_2}}\g  M_\varepsilon \right)\\
 	&=
		\frac{\int_{\Wt}e^{\frac{\tau}{\varepsilon}\langle f-f^\dagger,\psi\rangle_{W_1}
		+\tau\langle A\pt, \W\rangle_{W_2}+\ell(f_\tau)-\ell(f_\tau)
		+\ell(f)}d\Pi^{D_\varepsilon}(f)}{\int_{\Wt}e^{\ell(f)}d\Pi^{D_\varepsilon}(f)}. 
\end{align*}  
Using \eqref{Likelihood} we see that under $P^M_{f^\dagger}$
$$
	\ell(f)-\ell(f_\tau) 
	 =
		\frac{\tau}{\varepsilon}\langle A(f-f^\dagger),A\pt\rangle_{W_2}
		+\frac{\tau^2}{2}\|A\pt\|_{W_2}^2-\tau\langle A\pt,\W \rangle_{W_2}
$$
and hence 
\begin{align}\label{firstExp}
\E^{\Pi^{D_\varepsilon}}\left(e^{\frac{\tau}{\varepsilon}(\langle f,\psi\rangle_{W_1}-\widehat{\Psi})}\g M_\varepsilon \right) 
= e^{\frac{\tau^2}{2}\|A\pt\|_{W_2}^2}\frac{\int_{D_\varepsilon}e^{\ell(f_\tau)}d\Pi(f)}{\int_{D_\varepsilon}e^{\ell(f)}d\Pi(f)}.
\end{align}

	Let $\Pi_\tau$ be the shifted law of $f_\tau$, $f\sim\Pi$. Then by the Cameron-Martin theorem \cite[Corollary 2.4.3.]{Bogachev1998} we get, denoting $D_{\varepsilon,\tau}=\{g=f_\tau : f\in D_\varepsilon\}$,
\begin{align}
\label{eq:CMratio}
	\frac{\int_{D_{\varepsilon,\tau}}e^{\ell(g)}\frac{d\Pi_\tau(g)}{d\Pi(g)}d\Pi(g)}
	{\int_{D_\varepsilon}e^{\ell(g)}d\Pi(g)}
	=
		 \frac{\int_{D_{\varepsilon,\tau}}e^{\ell(g)}
		 e^{\tau\varepsilon\langle \pt,g\rangle_{V_\Pi}-\frac{(\tau\varepsilon)^2}{2}
		 \|\pt\|_{V_\Pi}^2}d\Pi(g)}{\int_{D_\varepsilon}e^{\ell(g)}d\Pi(g)}.
\end{align}
Since $\tilde\psi$ is a fixed element in $V_\Pi$, we see that $\varepsilon^2\|\pt\|_{V_\Pi}^2\to0$ as $\varepsilon\to0$. Using the definition of $D_\varepsilon$ in \eqref{simpleDeps} we see, as $\varepsilon\to0$,
\begin{align*}
	\varepsilon\sup_{g\in D_{\varepsilon,\tau}}|\langle \pt, g\rangle_{V_\Pi}|
	=
		\varepsilon\sup_{f\in D_{\varepsilon}}|\langle \pt, f+\tau\varepsilon\pt\rangle_{V_\Pi}|
	\leq
		 T\delta_\varepsilon\|\pt\|_{V_\Pi} +|\tau|\varepsilon^2\|\pt\|_{V_\Pi}^2
		 \to0.
\end{align*}
We have thus shown that a small shift of $f$ along $V_\Pi$ in \eqref{firstExp} correspond asymptotically to a shift in $D_\varepsilon$:
\begin{align*}
	\frac{\int_{D_\varepsilon}e^{\ell(f_\tau)}d\Pi(f)}{\int_{D_\varepsilon}e^{\ell(f)}d\Pi(f)}
	=
		\frac{\int_{D_{\varepsilon,\tau}}e^{\ell(g)}d\Pi(g)}{\int_{D_\varepsilon}
		e^{\ell(g)}d\Pi(g)}(1+o(1)) 
	= 
		\frac{\Pi(D_{\varepsilon,\tau} \g M_\varepsilon)}{\Pi(D_\varepsilon \g M_\varepsilon)}		(1+o(1)).
\end{align*}
Using Lemma \ref{lem:ae} we see that $\Pi(D_\varepsilon \g M_\varepsilon)\to1$ in $P_{f^\dagger}^M$-probability. We also note that 
\begin{align*}
	\Pi(D^c_{\varepsilon,\tau}) 
	& = 
		\Pi \left( g:  \frac{|\langle \pt, g-\tau\varepsilon\pt\rangle_{V_\Pi}|}{\|\pt\|_{V_\Pi}}
		 >\frac{T\delta_\varepsilon}{\varepsilon}\right)\\
	& \leq 
		\Pi \left( g: \frac{|\langle \pt, g\rangle_{V_\Pi}|}{\|\pt\|_{V_\Pi}}
		 >\frac{T\delta_\varepsilon}{\varepsilon}-|\tau|\varepsilon\|\pt\|_{V_\Pi}^2\right)\\
 	&\leq 
		e^{-\frac{t^2}{2}(\delta_\varepsilon/\varepsilon)^2}, 
\end{align*} 
for any $\sqrt{6}<t<T$. Using Lemma \ref{lem:ae} again we then conclude that $\Pi(D_{\varepsilon,\tau} \g M_\varepsilon)\to 1$ in $P_{f^\dagger}^M$-probability. 

%
%
%

\section*{3 $\ \ \ $ Step V in the proof of Theorem \ref{Mtheo}: convergence of the moments}
\label{MomentConv}

	The last step consists in replacing, in the result derived in the previous step, the centring $\hat\Psi_\varepsilon$ (defined in (\ref{hatpsieps})) with the posterior mean $\bar f=\E^{\Pi}[f|M_\varepsilon]$. The proof only requires minor adjustments from the proof of Theorem 2.7 in \cite{Monard2019}.  In particular, we show that 
\eq{
\label{postmeanconv}
	\|\varepsilon^{-1}(\bar f-\hat\Psi_\varepsilon)\|_{(H^\beta_K)^*}
	=
		\|\E^\Pi[\varepsilon^{-1}(f-\hat\Psi_\varepsilon)|M_\varepsilon]\|_{(H^\beta_K)^*}		=
			o_{P^M_{f^\dagger}}(1),
}
as $\varepsilon\to0$. 

	W argue by contradiction: let $(\Omega, \Sigma,\P)$ be the probability space on which 
$
M_\varepsilon
=L^{-1}f^\dagger+\varepsilon\W
$
is defined, and assume that for some $\Omega'\in\Sigma,\ \Pr(\Omega')>0,$ and $\xi>0$, we have along a certain vanishing sequence $(\varepsilon_n)_{n\ge1}$
\eq{
\label{nonconvmean}
	\|\E^\Pi[\varepsilon_n^{-1}(f-\hat\Psi_{\varepsilon_n})|
	M_{\varepsilon_n}(\omega)]\|_{(H^\beta_K)^*}\ge\xi, 
	\quad \forall\omega\in\Omega'.
}
In view of the convergence established in Step IV and since convergence in probability implies almost sure convergence for a subsequence, we can find $\Omega_0\in\Sigma,\  \Pr(\Omega_0)=1,$ such that along a further subsequence (denoted again as $(\varepsilon_n)_{n\ge1}$ for convenience)
$$
	\beta_{(H^\beta_K)^*}(\LC(\varepsilon_n^{-1}(f-\hat\Psi_{\varepsilon_n})|
	M_{\varepsilon_n}(\omega)),\mu)\to0
	\quad \textnormal{as}\ n\to\infty \ \forall \omega\in\Omega_0.
$$
Thus, for each $\omega\in\Omega_0$, recalling the definition (\ref{hatXeps}) of the process $\hat X_\varepsilon$ with law $\LC(\varepsilon^{-1}(f-\hat\Psi_\varepsilon)|M_\varepsilon)$ on $(H^\beta_K)^*$, the sequence $\{\hat X_{\varepsilon_n}(\omega), \ n\ge1\}$ of Borel random elements in $(H^\beta_K)^*$ will convergence in distribution to the process $X$  in (\ref{X}).
By Skorohod's embedding theorem \cite[Theorem 11.7.2]{Dudley2002} we can find a probability space and random elements with values in $(H^\beta_K)^*$, $\tilde X_{\varepsilon_n}(\omega)\stackrel{d}{=}\hat X_{\varepsilon_n}(\omega), \ \tilde X\stackrel{d}{=}X,$ defined on it such that
$
	\tilde X_{\varepsilon_n}(\omega)\to^{\textnormal{a.s.}}\tilde X,
$
or, equivalently,
\eq{\label{asconv}
	\|\tilde X_{\varepsilon_n}(\omega)-\tilde X\|_{(H^\beta_K)^*} \to^{\textnormal{a.s.}}0.
}

	From the standard conjugacy property of Gaussian priors with respect to linear inverse problems with Gaussian noise, $\hat X_{\varepsilon_n}(\omega)$ is a Gaussian random element in $(H^\beta_K)^*$  for each $\omega\in\Omega_0, \ n\ge1$. Then, also $\hat X_{\varepsilon_n}(\omega)-\tilde X$ is Gaussian, and by the Paley-Zygmund argument in Exercise 2.1.4 in \cite{Gine2016}, (\ref{asconv}) in fact implies the convergence of all norm-moments; in particular:
$$
	\E^\Pi\|\tilde X_{\varepsilon_n}(\omega)-\tilde X\|_{(H^\beta_K)^*}\to 0.
$$
Thus, since $X$ is a centred process, we obtain that for each $\omega\in\Omega_0$
\eqstar{
	\|\E^\Pi[\varepsilon_n^{-1}(f-\hat\Psi_{\varepsilon_n})|M(\omega)]\|_{(H^\beta_K)^*}
	&=
		 \|\E\tilde X_{\varepsilon_n}(\omega)-\E \tilde X \|_{(H^\beta_K)^*}
	\to0,
}
contradicting (\ref{nonconvmean}) since $\Pr(\Omega_0)=1$.

\end{document}